\documentclass[paper=A6,12pt]{article}

\usepackage{amsmath}
\usepackage{verbatim}
\usepackage{amssymb}
\usepackage{latexsym}
\usepackage{color}
\definecolor{green}{rgb}{0,1,0}
\usepackage{cancel}

\usepackage{xurl}
\usepackage{xcolor}
\definecolor{linkColor}{RGB}{16,16,16}
\usepackage{hyperref}
\hypersetup{colorlinks,allcolors=linkColor}

\hypersetup{
           breaklinks=true,   % splits links across lines
           colorlinks=true,   % displays links as colored text instead of blocks
           pdfusetitle=true,  % \title and \author values into pdf metadata
        }
\begin{document}
\newfont{\blb}{msbm10 scaled\magstep1} 

\newtheorem{theo}{Theorem}[section]

\newtheorem{defi}[theo]{Definition}
\newtheorem{prop}[theo]{Proposition}
\newtheorem{lemm}[theo]{Lemma}
\newtheorem{coro}[theo]{Corollary}
\pagestyle{myheadings}
%\markright{Torsion free groups with all subgroups $4$-subnormal}
\date{}
\author{Layla Sorkatti \\
Gunnar Traustason \\
Department of Mathematical Sciences \\ 
University of Bath, UK}

\title{Nilpotent symplectic alternating algebras II}
\maketitle
\begin{abstract}
\noindent
In this paper and its sequel we continue our study of nilpotent symplectic alternating
algebras. In particular we give a full classification of such algebras of dimension $10$ over any field. It is known that symplectic alternating algebras over $\mbox{GF}(3)$ correspond to a special rich class $\mathcal{C}$ of $2$-Engel $3$-groups of exponent $27$ and under this correspondence we
will see that the nilpotent algebras correspond to a subclass of $\mathcal{C}$ that are those groups in $\mathcal{C}$ that have an extra group theoretical property that we refer to as being powerfully
nilpotent and can be described also in the context of $p$-groups where $p$ is an arbitrary
prime.
\end{abstract}
\section{Introduction}
A symplectic alternating algebra (SAA) is a symplectic vector space L, whose associated alternating form is non-degenerate, that is furthermore equipped with a binary alternating product $ \cdot $ $ : L \times L  \mapsto L$ with the extra requirement that
$$(x \cdot y, z) = (y \cdot z , x)$$
for all $x,y,z \in L$. This condition can be expressed equivalently by saying that $(u \cdot x , v) = (u, v\cdot x)$ for all $u,v,x \in L$ or in other words that multiplication from the right is self-adjoint with respect to the alternating form.\\ \\
Symplectic alternating algebras originate from a study of powerful $2$-Engel groups~\cite{moravec,gt-2008} and there is in a 1-1 correspondence between a certain rich class of powerful $2$-Engel $3$-groups of exponent $27$ and SAAs over the field $\mbox{GF}(3)$. We will consider this in more details later in the introduction.\\ \\
Let $2n$ be a given even integer and $\mathbb{F}$ a fixed field. Let $V$ be the symplectic vector space over the field $\mathbb{F}$ with a non-degenerate alternating form. 
Fix some basis $u_1, u_2, \ldots, u_{2n}$ for $V$. An alternating product $\cdot$ that turns $V$ into a symplectic alternating algebra is uniquely determined by the values
$$ \mathcal{P} : \quad (u_i \cdot u_j, u_k), \quad 1 \leq i < j < k \leq 2n.$$
Let $L$ be the resulting symplectic alternating algebra. We refer to the data above as a presentation for $L$ with respect to the basis $u_1, \ldots, u_{2n}$.\\ \\
If $m(n)$ is the number of symplectic alternating algebras over a finite field $\mathbb{F}$ then $m(n) = |\mathbb{F} |^{\frac{4n^3}{3}+ O(n^2)}$ [$2$]. Because of the sheer growth, a general classification does not seem to be within reach although this has been done for small values of $n$. Thus it is not difficult to see that $m(0) = m(1) = 1$ 
$\mbox{and } m(2) = 2$. For higher dimensions the classification is already difficult. It is though known that when $\mathbb{F} = \mbox{GF}(3)$ we have $m(3)=31$~\cite{nsaa1}. 
Some general structure theory is developed in ~\cite{saa, to, nsaa1, sor-thesis}.
In particular there is dichotomy result that is an analog to a corresponding theorem for Lie algebras, namely that $L$ either contains a non-trivial abelian ideal or is a direct sum of simple symplectic alternating algebras~\cite{saa}. We also have that any symplectic algebra that is abelian-by-nilpotent must be nilpotent while this is not the case in general for solvable algebras~\cite{to}.\\ \\ 
Here we focus on the subclass of nilpotent symplectic alternating algebras (NSAAs) and this paper is a sequel to~\cite{nsaa1} %[2]
where a rich general structure theory was developed for NSAAs with a number of beautiful properties. Before discussing these we need to introduce some notation. Firstly we can always pick a basis $x_1, y_1, \ldots, x_n, y_n$ with the property that $(x_i, x_j) = (y_i, y_j) =0$ and $(x_i, y_j) = \delta_{ij}$ for $1 \leq i \leq j \leq k \leq n$.  We refer to a basis of this type as a standard basis. It turns out that for any nilpotent symplectic alternating algebra one can always choose a suitable standard basis such that the chain of subspaces
$$ 0 = I_0 < I_1 < \ldots < I_n < I_{n-1}^\perp < \cdots < I_0^\perp = L,$$
with $I_k = \mathbb{F} x_n + \cdots + \mathbb{F} x_{n-k+1}$ for $k > 0$, is a central chain of ideals. One can furthermore see from this that $x_iy_j = 0$ if $j \leq i$ and that $I^\perp_{n-1}$ is abelian. It follows that a number of the triple values $(uv,w)$ are trivial. Listing only the values that are possibly non-zero it suffices to consider
\[\mathcal{P}:\quad (x_i y_j, y_k)= \alpha_{ijk},\ \ (y_iy_j, y_k)=\beta_{ijk}\]
for some $\alpha_{ijk}, \beta_{ijk} \in \mathbb{F}$ where $1 \leq i < j < k\leq n$. Such a presentation is called a \emph{nilpotent presentation}.
Conversely any such presentation describes a nilpotent SAA. The algebras that are of maximal class turn out to have a rigid ideal structure. In particular when $2n \geq 10$ we can choose our chain of ideals above such that they are all characteristic and it turns out that $I_0, I_2, I_3, \ldots, I_{n-1}, I^\perp_{n-1}, I^\perp_{n-2}, \ldots, I^\perp_0$ are unique and equal to both the terms of the lower and upper central series (see~\cite{nsaa1} %[2] 
Theorem 3.1 and 3.2). The algebras of maximal class can be identified easily from their nilpotent presentations. In fact, if $\mathcal{P}$ is any nilpotent presentation of L with respect to a standard basis $\{ x_1, y_1, \ldots, x_n, y_n  \}$, and $2n \geq 8$, we have that L is of maximal class if and only if $x_i y_{i+1} \neq 0$ for all $i=2, \ldots, n-2$, and $x_1y_2, y_1y_2$ are linearly independent (see~\cite{nsaa1} %[2]
Theorem 3.4).\\ \\
From the general theory of nilpotent SAAs one can also determine their growth. Thus if $k(n)$ is the number of nilpotent SAAs of dimension $2n$ over a finite field ${\mathbb F}$ then $k(n)=|{\mathbb F}|^{n^3/3 + O(n^2)}$~\cite{nsaa1} %[2].
Again the growth is too large for a general classification to be feasable. The algebras of dimension $2n$ for $n \leq 4$ are classified in~\cite{nsaa1} %[2] 
over any field.
In this paper and its sequel we deal with the challenging classification of algebras of dimension $10$ over any field. As we will see the classification depends very much on the underlying field. It turns out that the classification of nilpotent symplectic alternating algebras of dimension 10 with a centre that is not isotropic can be easily reduced to the known classification of algebras of dimension $8$.
The main bulk of the work is thus about algebras with isotropic centre that must lie between $2$ and $5$. 
In this paper we consider the situation when the isotropic centre is of dimension $3$ or $5$ leaving the remaining cases to the sequel. As we said above, the classification depends on the underlying field.
We can read in particular from the classification that over a field that is algebraically closed there are $22$ NSAAs of dimension $10$. \\ \\
\section{The correspondence between SAAs and groups}
\noindent
Before starting the work on the classification we consider here in more detail the correspondence between SAAs and $2$-Engel groups, mentioned in the introduction. The study in $[4]$ reveals that there is one-one correspondence between symplectic alternating algebras over the field $\mbox{GF}(3)$ and a certain class $\mathcal{C}$ of powerful $2$-Engel $3$-groups of exponent $27$. These groups form a class  that consists of all powerful $2$-Engel $3$-groups $G$ with the following extra properties:\\ \\
(1) $G = \langle x, H \rangle$ where $H=\{ g \in G: g^9=1 \} $ and $Z(G)=\langle x \rangle$ with $O(x)=27$.\\ 
(2) $G$ is of rank $2r+1$ and has order $3^{3+4r}$. \\ \\
The associated symplectic alternating algebra $L(G)$ is constructed as follows. First we consider $L(G) = H/G^3$ as a vector space over GF($3$). To this we associate a bilinear alternating form (,) and an alternating binary multiplication as follows: for any $\bar{a} = a G^3, \bar{b} = b G^3$ and $\bar{c} = c G^3 $ in $L(G)$,
\[ [a,b]^3 = x^{9 (\bar{a}, \bar{b})} \]
\[ \bar{a} \cdot \bar{b} = \bar{c} \mbox{  where  } [a,b] Z(G) = c^3 Z(G).\]
One can show that these are well defined and turn $L(G)$ into a SAA. Furthermore $L(G) \cong L(K)$ if and only if $G \cong K$~\cite{gt-2008}
.\\ \\
In order to indentify the groups in $\mathcal{C}$ that correspond to the SAAs that are nilpotent, we introduce some new terms.\\ \\
{\bf Definition}. A finite p-group $G$ is \emph{powerfully nilpotent} if there exist an ascending chain
\[  \{1\} = H_0 \leq H_1 \leq \cdots \leq H_n=G \]
such that $[H_i, G] \leq H_{i-1}^p$ for $i=1, \cdots, n$. We refer to such a chain as a powerfully central chain and $n$ is the length of the chain. If $G$ is powerfully nilpotent then the smallest possible length of a powerfully central chain for $G$ is called its \emph{powerful nilpotence class}. \\ \\
Let us now turn to our special class $\mathcal{C}$ of powerful $2$-Engel $3$-groups. Let $G \in \mathcal{C}$. For any $K$ such that $G^3 \leq K \leq G$ we let $\bar{K} = K/G^3$. Notice that
\[ \bar{A} \cdot L(G) \leq \bar{B} \mbox{ if and only if } [ \langle A, x \rangle, G] \leq \langle B, x \rangle^3 .\]
Thus if $G^3 \leq H_i$ for $i=1, \cdots, n$, then
\[ \{0\} = \bar{H_0} \leq \bar{H_1} \leq \cdots \leq \bar{H_n}=L(G) \]
is a central chain of ideals in $L(G)$ if and only if 
\[ \{1\} \leq \langle x \rangle \leq \langle H_0, x \rangle \leq \cdots \leq \langle H_n, x \rangle = G \]
is a powerfully central chain. The classification of the NSAAs of dimension $10$ over $\mbox{GF}(3)$ gives us thus the classification for the powerfully nilpotent groups in $\mathcal{C}$ that are of rank $11$. The classification reveals that there are $25$ such groups.
\section{Algebras with a non-isotropic centre and algebras with an isotropic centre of dimension $5$}
\noindent
We consider first  the algebras with a non-isotropic centre. Let $L$ be such an algebra. In this case we can assume that we have a standard basis where $x_5,y_5 \in Z(L)$. We then have that $L$ is a direct sum of the  abelian algebra 
$\mathbb{F}x_5 + \mathbb{F}y_5$ 
and a NSAA of dimension $8$. The algebras of dimension $8$ were however classified in $[2]$ and according to this classification there are, apart from the abelian algebra, two algebras and one family of algebras. From this we can read that there are the following non-abelian NSAAs of dimension $10$ with a non-isotropic centre.
\begin{eqnarray*}
& {\mathcal Q}_{10}^{(7,1)}: \ \ 
 &(y_1y_2, y_3) = 1. \\
&{\mathcal Q}_{10}^{(5,1)}: \ \ 
&(y_1y_2, y_3) = 1,\ (x_1y_3, y_4) = 1. \\
&{\mathcal Q}_{10}^{(4,1)}(r): \ \ 
&(x_2 y_3, y_4) = r,\ (x_1y_2, y_4) = 1,\ (y_1y_2, y_3)=1,
\end{eqnarray*}
where $r \in \mathbb{F} \setminus \{ 0 \}$ and ${{\mathcal Q}_{10}}^{(4,1)}(s) \cong {{\mathcal Q}_{10}}^{(4,1)}(r)$ if and only if $r/s \in (\mathbb{F}^*)^3$. Here the notation ${{\mathcal Q}_{10}}^{(m,1)}$ indicates that the algebra has dimension $10$ with centre of dimension $m$. From now on we can thus assume that all our algebras have an isotropic centre and we start considering the case when the centre has dimension $5$.
Let $L$ be a nilpotent SAA of dimension $10$ with an isotropic centre of dimension $5$.
We can then choose a standard basis $x_{1},y_{1},x_{2},y_{2},x_{3},y_{3},
x_{4},y_{4},x_{5},y_{5}$ such that
  $$Z(L)={\mathbb F}x_{5}+{\mathbb F}x_{4}+{\mathbb F}x_{3}+{\mathbb F}x_{2}
    +{\mathbb F}x_{1}.$$
Here $x_{1},y_{1},x_{2},y_{2},x_{3},y_{3},x_{4},y_{4},x_{5},y_{5}$ will be 
determined later such that some further conditions hold. 
The elements $y_{1},\ldots ,y_{5}$ are not in $Z(L)$ and without loss of generality we can assume that $(y_{1}y_{2},y_{3})=1$. Now suppose that 
$(y_{i}y_{j},y_{4})=\alpha_{ij}$ and $(y_{i}y_{j},y_{5})=\beta_{ij}$ for $1\leq i,j\leq 3$. Replacing $x_{1},x_{2},x_{3},y_{4},y_{5}$ by 
$$\begin{array}{ll}
  \tilde{x}_{1}=x_{1}+\alpha_{23}x_{4}+\beta_{23}x_{5}, & \tilde{y}_{4}=y_{4}-\alpha_{12}y_{3}-
                           \alpha_{23}y_{1}-\alpha_{31}y_{2}, \\
  \tilde{x}_{2}=x_{2}+\alpha_{31}x_{4}+\beta_{31}x_{5}, & \tilde{y}_{5}=y_{5}-\beta_{12}y_{3}-
                           \beta_{23}y_{1}-\beta_{31}y_{2}, \\ 
  \tilde{x}_{3}=x_{3}+\alpha_{12}x_{4}+\beta_{12}x_{5}, & \mbox{} 
\end{array}$$            
we can assume that our standard basis has the further property that 
$(y_{i}y_{j},y_{4})=(y_{i}y_{j},y_{5})=0$ for $1\leq i<j\leq 3$. As $y_{4}\not
\in Z(L)$, we know that one of $(y_{1}y_{4},y_{5}), (y_{2}y_{4},y_{5}),
(y_{3}y_{4},y_{5})$ is nonzero. Without loss of generality we can assume that
$(y_{1}y_{4},y_{5})=1$. The only triples whose values are not known are then 
$\alpha=(y_{2}y_{4},y_{5})$ and $\beta=(y_{3}y_{4},y_{5})$. Replacing 
$x_{1},y_{2},y_{3}$ by $\tilde{x}_{1}=x_{1}+\alpha x_{2}+\beta x_{3}$, $\tilde{y}_{2}=y_{2}-\alpha y_{1},\tilde{y}_{3}
=y_{3}-\beta y_{1}$, we get a new
standard basis where the only nonzero triple values are $(y_{1}y_{2},y_{3})=1$
and $(y_{1}y_{4},y_{5})=1$. We have thus proved the following result. 
\begin{prop} There is a unique nilpotent SAA of dimension $10$ that has isotropic centre of
dimension $5$. This algebra can be described by the nilpotent presentation 
    $${\mathcal P}_{10}^{(5,1)}:\ (y_{1}y_{2},y_{3})=1,\ \ (y_{1}y_{4},y_{5})=1.$$
\end{prop}
\section{Algebras with an isotropic centre of dimension
$3$}  
In this section we will be assuming that $Z(L)$ is isotropic of dimension
$3$. First we derive some properties that hold for these algebras. Here 
throughout 
\begin{eqnarray*}
Z(L) & = & {\mathbb F}x_{5}+{\mathbb F}x_{4}+{\mathbb F}x_{3} \\
L^{2} & = & Z(L)+{\mathbb F}x_{2}+{\mathbb F}x_{1}+{\mathbb F}y_{1}+
{\mathbb F}y_{2}.
\end{eqnarray*}
\begin{lemm} $Z(L)\leq L^{3}$.
\end{lemm}
{\bf Proof}\ \ Otherwise $Z_{2}(L)=(L^{3})^{\perp}\nleq Z(L)^{\perp}=L^{2}$.
Without loss of generality we can suppose that $y_{3}\in Z_{2}(L)\setminus L^{2}$. As $Z_{2}(L)\cdot L^{2}=\{0\}$, we then have $y_{3}\cdot L^{2}=\{0\}$. Now also $x_{2}\cdot L^{2}=\{0\}$. Let $\alpha =(x_{2}y_{4},y_{5})$ and $\beta=(y_{3}y_{4},y_{5})$. Notice that
$\alpha,\beta \not =0$ as $x_{2},y_{3}\not\in Z(L)$. But then 
      $$((\beta x_{2}-\alpha y_{3})y_{4},y_{5})=0$$
that implies that $\beta x_{2}-\alpha y_{3}\in Z(L)$. This is absurd. $\Box$
\begin{lemm} $\mbox{dim\,}L^{3}\geq 5$. 
\end{lemm} 
{\bf Proof}\ \ Otherwise $\mbox{dim\,}L^{3}\leq 4$ and as $Z(L)\leq L^{3}\leq
L^{2}=Z(L)^{\perp}$ we can choose our standard basis such that 
$Z(L)={\mathbb F}x_{5}+{\mathbb F}x_{4}+{\mathbb F}x_{3}$ and 
       $$L^{3}\leq {\mathbb F}x_{5}+{\mathbb F}x_{4}+
  {\mathbb F}x_{3}+{\mathbb F}x_{2}.$$
This implies that ${\mathbb F}x_{5}+{\mathbb F}x_{4}+{\mathbb F}x_{3}+
{\mathbb F}x_{2}+{\mathbb F}x_{1}+{\mathbb F}y_{1}\leq 
(L^{3})^{\perp}=Z_{2}(L)$ and (notice that $Z_{2}(L)\leq L^{2}$ as
$Z(L)\leq L^{3}$) $L^{2}=Z_{2}(L)+Fy_{2}$ that implies that
$L^{2}$ is abelian. Then for any $x\in L^{2}$ and $a,b,c\in L$, we have
   $$(x,abc)=-(x(ab),c)=-(0,c)=0$$
and $L^{3}\leq (L^{2})^{\perp}=Z(L)$. Hence $L^{3}=Z(L)$ and $Z_{2}(L)=L^{2}$. 
Suppose $L=Z_{2}(L)+{\mathbb F}u_{1}+{\mathbb F}u_{2}+{\mathbb F}u_{3}$. Then 
$L^{2}=Z(L)+
{\mathbb F}u_{1}u_{2}+{\mathbb F}u_{1}u_{3}+{\mathbb F}u_{2}u_{3}$ and we get the 
contradiction that $4=\mbox{dim\,}L^{2}-
\mbox{dim\,}Z(L)\leq 3$. $\Box$
\begin{lemm} If $\mbox{dim\,}L^{3}=5$, then $L^{3}$ is isotropic.
\end{lemm}
{\bf Proof}\ \ Otherwise we can choose our basis such that
$L^{3}={\mathbb F}x_{5}+{\mathbb F}x_{4}+{\mathbb F}x_{3}+{\mathbb F}x_{1}+
{\mathbb F}y_{1}$ and then 
$Z_{2}(L)=(L^{3})^{\perp}={\mathbb F}x_{5}+{\mathbb F}x_{4}+{\mathbb F}x_{3}+
{\mathbb F}x_{2}+{\mathbb F}y_{2}$ and 
as $L^{2}\cdot Z_{2}(L)=\{0\}$, it follows that $x_{1}y_{2}=y_{1}y_{2}
=0$. Then $L^{2}$ is abelian and thus we get the contradiction that $L^{3}\leq Z(L)$. $\Box$
\begin{lemm} $Z(L)\leq L^{4}$. 
\end{lemm}
{\bf Proof}\ \ We have seen that $\mbox{dim\,}L^{3}\geq 5$. So we can
choose our standard nilpotent basis such that either 
          $$L^{3}={\mathbb F}x_{5}+{\mathbb F}x_{4}+{\mathbb F}x_{3}+
{\mathbb F}x_{2}+{\mathbb F}x_{1}$$
or 
         $$L^{3}={\mathbb F}x_{5}+{\mathbb F}x_{4}+{\mathbb F}x_{3}+
{\mathbb F}x_{2}+{\mathbb F}x_{1}+{\mathbb F}y_{1}.$$
We consider the two cases in turn beginning with the first case. If 
$Z(L) \not \leq L^{4}$, then $\mbox{dim\,}Z(L)\cap L^{4}\leq 2$
and thus $\mbox{dim\,}L^{2}+Z_{3}(L)=\mbox{dim\,}(Z(L)\cap L^{4})^{\perp}\geq 8$. Suppose $L=L^{2}+Z_{3}(L)+{\mathbb F}u+{\mathbb F}v$. Then 
$L^{2}=L^{3}+Z_{2}(L)+{\mathbb F}uv=L^{3}+{\mathbb F}uv$ and we get the contradiction that
$\mbox{dim\,}L^{2}\leq 5+1=6$. We now turn to the second case
where $L^{3}={\mathbb F}x_{5}+{\mathbb F}x_{4}+{\mathbb F}x_{3}+
{\mathbb F}x_{2}+{\mathbb F}x_{1}+{\mathbb F}y_{1}$. We argue by
contradiction and suppose that $Z(L)\cap L^{4}<Z(L)$. Then we can
choose our basis such that 
     $$L^{4}\leq {\mathbb F}x_{5}+{\mathbb F}x_{4}+{\mathbb F}x_{2}$$
and $Z(L)\cap L^{4}\leq {\mathbb F}x_{5}+{\mathbb F}x_{4}$. 
Now $y_{3}\in (L^4)^{\perp}=
Z_{3}(L)$ and as $Z_{3}(L)\cdot L^{3}=\{0\}$, it follows that 
     $$x_{1}y_{3}=x_{2}y_{3}=y_{1}y_{3}=0.$$
It follows from this that  $x_{1}y_{2},
y_{1}y_{2},y_{3}y_{2}\in {\mathbb F}x_{5}+{\mathbb F}x_{4}$. Thus in particular these three 
elements
are linearly dependent and we have $(\alpha x_{1}+\beta y_{1}+
\gamma y_{3})y_{2}=0$ where not all of $\alpha,\beta, \gamma$ are
zero. Then $x_{2},\alpha x_{1}+\beta y_{1}+\gamma y_{3}$ commute with
all the basis elements except possibly $y_{4}$ and $y_{5}$. Suppose
\begin{eqnarray*}
      (x_{2}y_{4},y_{5}) & = & r \\
   ((\alpha x_{1}+\beta y_{1}+\gamma y_{3})y_{4},y_{5}) & = & s.
\end{eqnarray*}
If $r=0$ then we get the contradiction that $x_{2}\in Z(L)$ and if
$r\not =0$, we get the contradiction that $-sx_{2}+r\alpha x_{1}+
r\beta y_{1}+r\gamma y_{3}\in Z(L)$. $\Box$
\mbox{}\\ \\ 
After these more general results we classify all the algebras where
$Z(L)$ is isotropic of dimension $3$. We consider the two subcases 
$\mbox{dim\,}L^{3}=5$ and $\mbox{dim\,}L^{3}=6$ separately.
\subsection{The algebras where $\mbox{dim\,}L^{3}=5$}
We have seen that $L^{3}$ must be isotropic and 
thus in particular we have that $L^{3}=(L^{3})^{\perp}=Z_{2}(L)$ that implies
that  $L^{4}\leq Z(L)$. By Lemma 4.4 we thus have $L^{4}=Z(L)$.
We have thus determined the terms of the lower and the upper central
series
$$\begin{array}{ll}
\begin{tabular}{c|c|c|c} 
\cline{2-3} 
  $L^{2}\cdot L^{2}$ & $x_{5}$ & $y_{5}$ & \mbox{} \\
\cline{2-3}
   \mbox{} & $x_{4}$ & $y_{4}$ & $(L^{2}\cdot L^{2})^{\perp}$ \\
\mbox{}$Z(L)=L^{4}$  & $x_{3}$ & $y_{3}$ & \mbox{} \\ 
\cline{2-3} 
 \mbox{} & $x_{2}$ & $y_{2}$ & $Z_{3}(L)=L^{2}$ \\
\mbox{}$Z_{2}(L)=L^{3}$  \mbox{} & $x_{1}$ & $y_{1}$ & \mbox{} \\
\cline{2-3}
\end{tabular} &  
$$\begin{array}{l}
 Z(L)=L^{4}={\mathbb F}x_{5}+{\mathbb F}x_{4}+{\mathbb F}x_{3} \\ 
 Z_{2}(L)=L^{3}=Z(L)+{\mathbb F}x_{2}+{\mathbb F}x_{1} \\
 Z_{3}(L)=L^{2}=Z_{2}(L)+{\mathbb F}y_{1}+{\mathbb F}y_{2} 
\end{array}$$
\end{array}$$
{\bf Remark}. As $L^{2}\cdot Z_{2}(L)=\{0\}$ we see that $x_{1}y_{2}=0$. 
Now $L^{2}$ is not abelian as this would imply that $L^{3}\leq Z(L)$. It follows
that $y_{1}y_{2}\not =0$ and we get a
one dimensional characteristic subspace
       $$L^{2}\cdot L^{2}={\mathbb F}y_{1}y_{2}.$$
Notice that $y_{1}y_{2}\in Z(L)$. We choose our basis such that $y_{1}y_{2}=
x_{5}$. We will also work with the $9$ dimensional characteristic subspace 
    $$V=(L^{2}\cdot L^{2})^{\perp}={\mathbb F}x_{5}+{\mathbb F}x_{4}+
{\mathbb F}x_{3}+{\mathbb F}x_{2}+{\mathbb F}x_{1}+{\mathbb F}y_{1}+
{\mathbb F}y_{2}+{\mathbb F}y_{3}+{\mathbb F}y_{4}.$$
As $x_{1}y_{2}=0$ we have that  $x_{1}y_{3},x_{1}y_{4}\perp y_{2}$. As 
$y_{1}y_{2}=x_{5}$ we also have that $y_{1}y_{3},y_{1}y_{4}\perp y_{1},y_{2}$ and
$y_{2}y_{3},y_{2}y_{4}\perp y_{1},y_{2}$. It follows that
    $$V^{2}+L^{4}=(L^{2}+{\mathbb F} y_{3}+{\mathbb F} y_{4})(L^{2}+{\mathbb F}y_{3}+{\mathbb F}y_{4})+
L^{4}=
        {\mathbb F}y_{3}y_{4}+L^{4}.$$
We consider few subcases. 
\subsubsection{Algebras where $V^{2}\leq L^{4}$}
Notice then that $y_{3}y_{4}\in {\mathbb F}x_{5}$ and thus $x_{1}y_{3},x_{2}y_{3},x_{1}y_{4},
x_{2}y_{4}\in {\mathbb F}x_{5}$. As 
$L^{3}=L^{4}+{\mathbb F}x_{2}+{\mathbb F}x_{1}$, we have 
   $$L^{4}=({\mathbb F}x_{2}+{\mathbb F}x_{1})({\mathbb F}y_{3}+{\mathbb F}y_{4}+
{\mathbb F}y_{5})={\mathbb F}x_{5}+{\mathbb F}x_{2}y_{5}+
{\mathbb F}x_{1}y_{5}.$$
Pick $x_{5},y_{5},x_{2},x_{1},y_{1}$ satisfying the conditions above and
let 
    $$x_{4}=-x_{2}y_{5},\ x_{3}=-x_{1}y_{5}$$
We can then extend  $x_{5},x_{4},x_{3},y_{1},$
$y_{2},y_{5}$ 
to a standard basis $x_{5},x_{4},
x_{3},x_{2},x_{1},y_{1},$
$y_{2},y_{3},y_{4},y_{5}$ satisfying the conditions
above. All triples involving
both $x_{1}$ and $y_{2}$ are $0$. The remaining ones are
$$\begin{array}{lll}
  (x_{1}y_{3},y_{5})=1 & (x_{1}y_{3},y_{4})=0 & (x_{1}y_{4},y_{5})=0 \\
 (x_{2}y_{3},y_{4})=0 & (x_{2}y_{3},y_{5})=0 & (x_{2}y_{4},y_{5})=1 \\
 (y_{1}y_{2},y_{3})=0 & (y_{1}y_{2},y_{4})=0 & (y_{1}y_{2},y_{5})=1 \\
 (y_{2}y_{3},y_{4})=0 & (y_{2}y_{3},y_{5})=\alpha & (y_{2}y_{4},y_{5})=\beta \\
(y_{1}y_{3},y_{4})=0 & (y_{1}y_{3},y_{5})=\gamma & (y_{1}y_{4},y_{5})=\delta \\
(y_{3}y_{4},y_{5})=r & \mbox{} & \mbox{} 
\end{array}$$
Now let 
$$\begin{array}{lll}
 \tilde{y}_{3}  =  y_{3}+\alpha y_{1}-\gamma y_{2}-sx_{2}-s\gamma x_{3}-
s\delta x_{4}  && \tilde{y}_{2}  =  y_{2}-sx_{3}   \\
\tilde{x}_{1}  =  x_{1}-\alpha x_{3}-\beta x_{4} && 
\tilde{y}_{4}  = y_{4}+\beta y_{1} - \delta y_{2} \\
\tilde{x}_{2}  =  x_{2}+\gamma x_{3} + \delta x_{4}  &&
\end{array}$$ \noindent
\noindent
where $s=r+\alpha\delta-\beta\gamma$. 
One checks readily that we get a new standard basis with a presentation
like the one above where $\tilde{\alpha}=\tilde{\beta}=\tilde{\gamma}=
\tilde{\delta}=\tilde{r}=0$. \\ \\ 
So we arrive at a unique algebra with presentation \\ \\
${\mathcal P}_{10}^{(3,1)}:\ \ (x_{1}y_{3},y_{5})=1,\ (x_{2}y_{4},y_{5})=1,\ 
(y_{1}y_{2},y_{5})=1.$ \\ \\
One can check that the centre has dimension $3$ and that $L^{3}$ has dimension $5$. Also $((L^{2}\cdot L^{2})^{\perp})^{2}={\mathbb F} x_{5}\leq L^{4}$. 
\subsubsection{Algebras where $V^{2}\nleq L^{4}$ but $V^{2}\leq L^{3}$} 
Here we can pick our basis such that 
        $$V^{2}+L^{4}={\mathbb F}x_{5}+{\mathbb F}x_{4}+
{\mathbb F}x_{3}+{\mathbb F}x_{2}.$$
Notice that $V^{3}={\mathbb F}x_{2}y_{3}+{\mathbb F}x_{2}y_{4}$ and as $(y_{3}y_{4},x_{2})=0$
we have that $V^{3}\leq {\mathbb F}x_{5}$. As $x_{2}\not\in Z(L)$, we furthermore must
have that $\mbox{dim\,}V^{3}=1$. This means that there is a characteristic
ideal $W$ of codimension $1$ in $V$ such that $x_{2}W=V^{2}W=\{0\}$. We choose our
basis such that 
     $$W={\mathbb F}x_{5}+{\mathbb F}x_{4}+{\mathbb F}x_{3}+{\mathbb F}x_{2}+
{\mathbb F}x_{1}+{\mathbb F}y_{1}+{\mathbb F}y_{2}+{\mathbb F}y_{3}.$$
It follows that we have a chain of characteristic ideals:
$$\begin{array}{ll}
\begin{tabular}{c|c|c|c} 
\cline{2-3} 
  $L^{2}\cdot L^{2}$ & $x_{5}$ & $y_{5}$ & \mbox{} \\
\cline{2-3}
   $W^{\perp}$ & $x_{4}$ & $y_{4}$ & $V$ \\
\mbox{}$Z(L)=L^{4}$  & $x_{3}$ & $y_{3}$ & $W$ \\ 
\cline{2-3} 
 \mbox{}$V^{2}+L^{4}$ & $x_{2}$ & $y_{2}$ & $Z_{3}(L)=L^{2}$ \\
\mbox{}$Z_{2}(L)=L^{3}$  \mbox{} & $x_{1}$ & $y_{1}$ & $(V^{2})^{\perp}\cap L^{2}$ \\
\cline{2-3}
\end{tabular} &  
\begin{array}{l}
  L^{2}\cdot L^{2}={\mathbb F}x_{5} \\
   W^{\perp}={\mathbb F}x_{5}+{\mathbb F}x_{4} \\
 Z(L)=L^{4}={\mathbb F}x_{5}+{\mathbb F}x_{4}+{\mathbb F}x_{3} \\ 
  V^{2}+L^{4}={\mathbb F}x_{5}+{\mathbb F}x_{4}+{\mathbb F}x_{3}+
{\mathbb F}x_{2} \\
 Z_{2}(L)=L^{3}=Z(L)+{\mathbb F}x_{2}+{\mathbb F}x_{1} \\
 (V^{2})^{\perp}\cap L^{2}=L^{3}+{\mathbb F}y_{1} \\
 Z_{3}(L)=L^{2}=Z_{2}(L)+{\mathbb F}y_{1}+{\mathbb F}y_{2} \\
  W=L^{3}+{\mathbb F}y_{1}+{\mathbb F}y_{2}+{\mathbb F}y_{3} \\
  V=L^{3}+{\mathbb F}y_{1}+{\mathbb F}y_{2}+{\mathbb F}y_{3}+
         {\mathbb F}y_{4}
\end{array}
\end{array}$$
We want to show that there is again a unique algebra satisfying these
conditions. We modify the basis and reach a unique presentation. Notice that 
$V^{2}W=\{0\}$ and $L^{2}\cdot Z_{2}(L)=\{0\}$ imply that 
      $$x_{1}y_{2}=x_{2}y_{3}=0.$$
We have also chosen our basis such that 
\begin{equation}
    y_{1}y_{2}=x_{5}.
\end{equation}  
Notice next that $x_{2}y_{3}=0$ implies that $x_{2}y_{4}$ is orthogonal
to $y_{3}$ and $y_{4}$ and thus $x_{2}y_{4}=ry_{5}$ where $r$ must
be nonzero as $x_{2}\not\in Z(L)$. By replacing $y_{4}$ and $x_{4}$
by $ry_{4}$ and $\frac{1}{r}x_{4}$, we can assume that
\begin{equation}
   x_{2}y_{4}=x_{5}.
\end{equation}
As $y_{3}y_{4}\in V^{2}\leq L^{3}$ and as $x_{1}y_{2}=0$ we have that $x_{1}y_{4}$ is orthogonal to $y_{2},y_{3},y_{4}$. Thus
$x_{1}y_{4}=\alpha x_{5}$ for some $\alpha \in {\mathbb F}$. Replacing $x_{1}, y_{2}$ by
$x_{1}-\alpha x_{2}$ and $y_{2}+\alpha y_{1}$ we get a new standard
basis where 
\begin{equation}
   x_{1}y_{4}=0.
\end{equation}
Notice that the change in $y_{2}$ does not affect (1). We next turn our attention to $x_{1}y_{3}$. As $x_{1}y_{2}=0$ and $x_{1}y_{4}=0$, we have that $x_{1}y_{3}$ is orthogonal to $y_{1},y_{2},y_{3},y_{4}$
and thus $x_{1}y_{3}=rx_{5}$ where $r$ is nonzero since $x_{1}\not\in Z(L)$. By replacing $y_{3}$ and $x_{3}$ by $ry_{3}$ and $\frac{1}{r}x_{3}$ we get 
\begin{equation}
       x_{1}y_{3}=x_{5}.
\end{equation}
Now we see, as $y_{1}y_{2}=x_{5}$ and $y_{3}y_{4}\in L^{4}+V^{2}$, that 
$y_{1}y_{3}$ is orthogonal to $y_{1},y_{2},y_{3}$ and $y_{4}$. Thus
$y_{1}y_{3}=ax_{5}$ for some $a\in  {\mathbb F}$. Replacing $y_{1}$ by $y_{1}-
ax_{1}$ we can assume that 
\begin{equation}
     y_{1}y_{3}=0.
\end{equation}
As $x_{1}y_{2}=0$ the change in $y_{1}$ does not affect (1). From the discussion above we know that $y_{1}y_{4}$ is orthogonal to $y_{1},y_{2},y_{3}$ and
$y_{4}$ and thus $y_{1}y_{4}=ax_{5}$ for some $a\in {\mathbb F}$. Replacing
$y_{4},x_{2}$ by $y_{4}-ay_{2}$ and $x_{2}+ax_{4}$, we get a new
standard basis where 
\begin{equation}
   y_{1}y_{4}=0.
\end{equation}
These changes do not affect (2) and (3). As $y_{3}y_{4}\in V^{2}+L^{4}$
but not in $L^{4}$ we know that $(y_{2}y_{3},y_{4})=r$ for some nonzero 
$r\in {\mathbb F}$. Suppose also that $(y_{3}y_{4},y_{5})=\alpha$. Then 
$y_{3}y_{4}=rx_{2}+\alpha x_{5}$. Replace $x_{2}$ and $y_{5}$
by $x_{2}+\frac{\alpha}{r}x_{5}$ and $y_{5}-\frac{\alpha}{r}y_{2}$. Then
\begin{equation}
   y_{3}y_{4}=rx_{2}.
\end{equation}
The changes do not affect (2). Then consider the triples 
    $$(y_{2}y_{3},y_{5})=a,\ \ (y_{2}y_{4},y_{5})=b.$$
Replacing $y_{5},x_{4},x_{3}$ by $y_{5}-\frac{a}{r}y_{4}+
\frac{b}{r}y_{3}$, $x_{4}+\frac{a}{r}x_{5}$ and $x_{3}-\frac{b}{r}x_{5}$ we can assume that 
\begin{equation}
   (y_{2}y_{3},y_{5})=(y_{2}y_{4},y_{5})=0.
\end{equation}
We have then arrived at a presentation where the only nonzero triples
are 
   $$(x_{2}y_{4},y_{5})=1,\ (x_{1}y_{3},y_{5})=1,\ (y_{1}y_{2},y_{5})=1,\ (y_{2}y_{3},y_{4})=r.$$
Replacing $x_{1},y_{1},x_{2},y_{2},x_{3},y_{3},x_{4},y_{4}$, by
$\frac{1}{r}x_{1},ry_{1},rx_{2},\frac{1}{r}y_{2},\frac{1}{r}x_{3},
ry_{3},rx_{4},\frac{1}{r}y_{4}$, we get a unique algebra with presentation: \\ \\
${\mathcal P}_{10}^{(3,2)}:\ (x_{2}y_{4},y_{5})=1,\ (x_{1}y_{3},y_{5})=1,\ 
(y_{1}y_{2},y_{5})=1,\ (y_{2}y_{3},y_{4})=1.$ \\ \\
One can easily check that conversely this algebra belongs to the category that we have been studying. 
\subsubsection{Algebras where $V^{2}\leq L^{2}$ but $V^{2}\nleq L^{3}$}
Pick our basis such that 
     $$V^{2}+L^{4}={\mathbb F} x_{5}+{\mathbb F} x_{4}+{\mathbb F} x_{3}+{\mathbb F} y_{1}$$
Notice then that 
    $${\mathbb F} x_{5}+{\mathbb F} x_{4}+{\mathbb F} x_{3}+{\mathbb F}x_{2}=(V^{2}+L^{4})^{\perp} \cap L^{3}=
    (V^{2})^{\perp}\cap L^{3}$$
is a characteristic ideal of $L$. As $y_{3}y_{4}\in V^{2}+L^{4}$ we have
that $x_{2}y_{3}\perp y_{4}$ and $x_{2}y_{4}\perp y_{3}$. Thus 
$x_{2}V\leq {\mathbb F}x_{5}$. As $x_{2}\not\in Z(L)$ , we must furthermore
have that $x_{2}V=((V^{2})^{\perp}\cap L^{3})V={\mathbb F}x_{5}$. This implies that the centraliser of $(V^{2})^{\perp}\cap L^{3}$ in $V$ is a characteristic ideal $W$
of codimension $1$. We can choose our basis such that
    $$W={\mathbb F} x_{5}+{\mathbb F} x_{4}+{\mathbb F} x_{3}+{\mathbb F} x_{2}+{\mathbb F} x_{1}+{\mathbb F} y_{1}+{\mathbb F} y_{2}+{\mathbb F}y_{3}.$$
We now get a chain of characteristic ideals as before
$$\begin{array}{ll}
\begin{tabular}{c|c|c|c} 
\cline{2-3} 
  $L^{2}\cdot L^{2}$ & $x_{5}$ & $y_{5}$ & \mbox{} \\
\cline{2-3}
   $W^{\perp}$ & $x_{4}$ & $y_{4}$ & $V$ \\
\mbox{}$Z(L)=L^{4}$  & $x_{3}$ & $y_{3}$ & $W$ \\ 
\cline{2-3} 
 \mbox{}$(V^{2})^{\perp}\cap L^{3}$ & $x_{2}$ & $y_{2}$ & $Z_{3}(L)=L^{2}$ \\
\mbox{}$Z_{2}(L)=L^{3}$  \mbox{} & $x_{1}$ & $y_{1}$ & $V^{2}+L^{3}$ \\
\cline{2-3}
\end{tabular} &  
\begin{array}{l}
  L^{2}\cdot L^{2}={\mathbb F}x_{5} \\
   W^{\perp}={\mathbb F}x_{5}+{\mathbb F}x_{4} \\
 Z(L)=L^{4}={\mathbb F}x_{5}+{\mathbb F}x_{4}+{\mathbb F}x_{3} \\ 
  (V^{2})^{\perp}\cap L^{3}={\mathbb F}x_{5}+{\mathbb F}x_{4}+{\mathbb F}x_{3}+
{\mathbb F}x_{2} \\
 Z_{2}(L)=L^{3}=Z(L)+{\mathbb F}x_{2}+{\mathbb F}x_{1} \\
 V^{2}+ L^{3}=L^{3}+{\mathbb F}y_{1} \\
 Z_{3}(L)=L^{2}=Z_{2}(L)+{\mathbb F}y_{1}+{\mathbb F}y_{2} \\
  W=L^{3}+{\mathbb F}y_{1}+{\mathbb F}y_{2}+{\mathbb F}y_{3} \\
  V=L^{3}+{\mathbb F}y_{1}+{\mathbb F}y_{2}+{\mathbb F}y_{3}+
         {\mathbb F}y_{4}
\end{array}
\end{array}$$
As $((V^{2})^{\perp}\cap L^{3})W=\{0\}$ and $L^{2}\cdot Z_{2}(L)=\{0\}$, we see that 
            $$x_{1}y_{2}=x_{2}y_{3}=0.$$
We have also chosen our basis such that 
\begin{equation}
             y_{1}y_{2}=x_{5}.
\end{equation}
Notice next that $x_{2}y_{3}=0$ implies that $x_{2}y_{4}$ is orthogonal to
$y_{3}$ and  $y_{4}$ and thus $x_{2}y_{4}=ry_{5}$ where $r$ must be nonzero as $x_{2}\not\in Z(L)$. By replacing $y_{4}$ and $x_{4}$ by $ry_{4}$ and $\frac{1}{r}x_{4}$, we can assume that 
\begin{equation}
   x_{2}y_{4}=x_{5}.
\end{equation}
As $y_{3}y_{4}\in V^{2}+L^{4}$ and $y_{1}y_{2}=x_{5}$, we have that $y_{1}y_{4}$ is orthogonal
to $y_{2},y_{3},y_{4}$. Thus $y_{1}y_{4}=ay_{5}$ for some $a\in {\mathbb F}$. Replacing 
$y_{4},x_{2}$ by $y_{4}-ay_{2}$ and $x_{2}+ax_{4}$ we get
\begin{equation}
   y_{1}y_{4}=0.
\end{equation}
Notice that the change does not affect (10). Next notice similarly that
$y_{1}y_{3}$ is orthogonal to $y_{2},y_{3},y_{4}$ and thus $y_{1}y_{3}=ay_{5}$ for some $a\in {\mathbb F}$. Replacing $y_{3}$ and $x_{2}$ by
$y_{3}-ay_{2}$ and $x_{2}+ax_{3}$ we get
\begin{equation}
   y_{1}y_{3}=0.
\end{equation}
Notice that (11) is not affected by this change. We know that 
$x_{1}y_{2}=0$. The possible nonzero triples involving $x_{1}$ are then 
 $$(x_{1}y_{3},y_{4})=r,\ (x_{1}y_{3},y_{5})=a,\ (x_{1}y_{4},y_{5})=b.$$
Notice that as $y_{3}y_{4}\in (Z(L)+{\mathbb F} y_{1})\setminus Z(L)$ we must have
that $r\not =0$. Replace $y_{5},x_{4},x_{3}$ by $y_{5}-\frac{a}{r}y_{4}
+\frac{b}{r}y_{3}$, $x_{4}+\frac{a}{r}x_{5}$ and $x_{3}-\frac{b}{r}x_{5}$ and we get a new standard basis where 
    $$x_{1}y_{3}=rx_{4},\ \ x_{1}y_{4}=-rx_{3}.$$
Replacing $y_{3},x_{3}$ by $ry_{3},\frac{1}{r}x_{3}$ gives 
\begin{equation}
    x_{1}y_{3}=x_{4},\ \ x_{1}y_{4}=-x_{3}.
\end{equation}
It follows that $(y_{2}y_{3},y_{4})=(y_{3}y_{4},y_{5})=0$. Suppose 
$(y_{2}y_{3},y_{5})=a,\ (y_{2}y_{4},y_{5})=b$. Replace $y_{3},y_{4},
x_{1}$ by $y_{3}+ay_{1},\ y_{4}+by_{1}$ and $x_{1}-ax_{3}-bx_{4}$. 
Notice that these changes do not affect the equations above and we
now arrive at a unique algebra with presentation: \\ \\
${\mathcal P}_{10}^{(3,3)}:\ (x_{2}y_{4},y_{5})=1,\ (x_{1}y_{3},y_{4})=1,\
(y_{1}y_{2},y_{5})=1.$ \\ \\
Calculations show that conversely this algebra belongs to the relevant
category. There are thus exactly three algebras where $Z(L)$ is isotropic of dimension $3$ and where the dimension of $L^{3}$ is $5$. 
\begin{prop} There are exactly three NSAAs of dimension $10$ that have an 
isotropic centre of dimension $3$ and where $\mbox{dim\,}L^{3}=5$. These
are given by the presentations: \\ \\
${\mathcal P}_{10}^{(3,1)}:\ \ (x_{1}y_{3},y_{5})=1,\ (x_{2}y_{4},y_{5})=1,\ 
(y_{1}y_{2},y_{5})=1.$ \\ \\
${\mathcal P}_{10}^{(3,2)}:\ (x_{2}y_{4},y_{5})=1,\ (x_{1}y_{3},y_{5})=1,\ 
(y_{1}y_{2},y_{5})=1,\ (y_{2}y_{3},y_{4})=1.$ \\ \\
${\mathcal P}_{10}^{(3,3)}:\ (x_{2}y_{4},y_{5})=1,\ (x_{1}y_{3},y_{4})=1,\
(y_{1}y_{2},y_{5})=1.$ 
\end{prop}
\subsection{The algebras where $\mbox{dim\,}L^{3}=6$} 
Here we are thus assuming that 
      $$L^{3}={\mathbb F}x_{5}+{\mathbb F}x_{4}+{\mathbb F}x_{3}+
   {\mathbb F}x_{2}+{\mathbb F}x_{1}+{\mathbb F}y_{1}.$$
\begin{lemm} We have $\mbox{dim\,}L^{4}=4$. 
\end{lemm}
{\bf Proof}. By Lemma 4.4, we know that $Z(L)\leq L^{4}$ and we also
$L^{4}\leq {\mathbb F}x_{5}+{\mathbb F}x_{4}+
{\mathbb F}x_{3}+{\mathbb F}x_{2}$. Thus if  the dimension of $L^{4}$ is not 
$4$, then $L^{4}=Z(L)={\mathbb F}x_{5}+{\mathbb F}x_{4}+{\mathbb F}x_{3}$ and 
$Z_{3}(L)=
(L^{4})^{\perp}=L^{3}+{\mathbb F}y_{2}$. As $Z_{3}(L)\cdot L^{3}=\{0\}$, it follows
that $x_{1}y_{2}=y_{1}y_{2}=0$ and $L^{2}$ is abelian. Hence we
get the contradiction that $L^{3}\leq Z(L)$. $\Box$ \\ \\
It follows that we have $L^{4}={\mathbb F}x_{5}+{\mathbb F}x_{4}+
{\mathbb F}x_{3}+{\mathbb F}x_{2}$. 
\begin{lemm}
We have $\mbox{dim\,}L^{5}=2$. 
\end{lemm}
{\bf Proof}\ \ We have an alternating form
    $$\phi:L/L^{2}\times L/L^{2}\rightarrow {\mathbb F}$$
given by $\phi(\bar{y},\bar{z})=(x_{2}y,z)$. As $L/L^{2}$ has odd dimension
we know that the isotropic part must be non-trivial. Thus we can then
choose our standard basis such that $(x_{2}y_{3},y_{4})=(x_{2}y_{3},y_{5})=0$
and thus $x_{2}y_{3}=0$. It follows that $L^{5}={\mathbb F}x_{2}(y_{3}+y_{4}+y_{5})={\mathbb F}x_{2}y_{4}+{\mathbb F}x_{2}y_{5}$ and thus it is of dimension at most $2$.
As $L^{4}\not\leq Z(L)$ we have $\mbox{dim\,}L^{5}>0$ and as we know [3, proposition 3.10] that
$\mbox{dim\,}L^{5}\not =1$ we must that $\mbox{dim\,}L^{5}=2$. $\Box$. \\ \\
We thus have determined the lower and upper central series of $L$. We have
$$\begin{array}{ll}
\begin{tabular}{c|c|c|c} 
\cline{2-3} 
 \mbox{} & $x_{5}$ & $y_{5}$ & \mbox{} \\
   $L^{5}$ & $x_{4}$ & $x_{4}$ & \mbox{} \\
\cline{2-3}
\mbox{}$Z(L)$  & $x_{3}$ & $y_{3}$ & $Z_{4}(L)$ \\ 
\cline{2-3} 
 \mbox{}$L^{4}=Z_{2}(L)$ & $x_{2}$ & $y_{2}$ & $L^{2}$ \\
\cline{2-3}
\mbox{} & $x_{1}$ & $y_{1}$ & $L^{3}=Z_{3}(L)$ \\
\cline{2-3}
\end{tabular} &  
$$\begin{array}{l}
   L^{5}={\mathbb F}x_{5}+{\mathbb F}x_{4} \\
 Z(L)={\mathbb F}x_{5}+{\mathbb F}x_{4}+{\mathbb F}x_{3} \\ 
 Z_{2}(L)=L^{4}={\mathbb F}x_{5}+{\mathbb F}x_{4}+{\mathbb F}x_{3}+
{\mathbb F}x_{2} \\
 Z_{3}(L)=L^{3}=Z(L)+{\mathbb F}x_{2}+{\mathbb F}x_{1}+{\mathbb F}y_{1} \\
 L^{2}=L^{3}+{\mathbb F}y_{2} \\
  Z_{4}(L)=L^{3}+{\mathbb F}y_{2}+{\mathbb F}y_{3}
\end{array}$$
\end{array}$$
Notice that $x_{2}\in L^{4}$ and $y_{3}\in Z_{4}(L)$ and thus $x_{2}y_{3}=0$.
Also 
     $$L^{3}L^{2}={\mathbb F}x_{1}y_{2}+{\mathbb F}y_{1}y_{2}\leq Z(L).$$
Furthermore $x_{1}y_{2}$ and $y_{1}y_{2}$ are linearly independent. To see this we argue by contradiction and suppose that $0=ax_{1}y_{2}+by_{1}y_{2}$ 
for some $a,b\in {\mathbb F}$ where not both $a,b$ are zero. Then $(ax_{1}+by_{1})L\leq Z(L)$ that would give us the contradiction that 
$ax_{1}+by_{1}\in Z_{2}(L)$. \\ \\
We thus have that $L^{3}L^{2}$ is a $2$-dimensional subspace of $Z(L)$ and 
we consider two possible cases namely $L^{3}L^{2}=L^{5}$ and 
$L^{3}L^{2}\not =L^{5}$. We consider the latter first.
\subsubsection{Algebras where $L^{3}\cdot L^{2}\neq L^{5}$}
Here $L^{3}L^{2}\cap L^{5}$ is one dimensional and we can choose our
standard basis such that $L^{3}L^{2}=L^{3}y_{2}={\mathbb F}x_{5}$. In order to 
clarify
the structure further we introduce the following isotropic characteristic ideal of dimension $5$:
      $$U=\{x\in L^{3}:\,xL^{2}\leq L^{3}L^{2}\cap L^{5}\}.$$
Now $L^{3}L^{2}$ is of dimension $2$ and $L^{4}L^{2}=0$ and thus $U$ is of codimension $1$ in $L^{3}$ and contains $L^{4}$. We can thus choose our standard
basis such that $U={\mathbb F}x_{5}+{\mathbb F}x_{4}+{\mathbb F}x_{3}+
{\mathbb F}x_{2}+{\mathbb F}x_{1}$. We thus have the following picture \\ \\
\begin{tabular}{c|c|c|c} 
\cline{2-3} 
 $L^{5}\cap L^{3}L^{2}$ & $x_{5}$ & $y_{5}$ & \mbox{} \\
\cline{2-3}
   $L^{5}$ & $x_{4}$ & $x_{4}$ & $(L^{3}L^{2}\cap L^{5})^{\perp}$ \\
\cline{2-3}
\mbox{}$Z(L)$  & $x_{3}$ & $y_{3}$ & $Z_{4}(L)$ \\ 
\cline{2-3} 
 \mbox{}$L^{4}=Z_{2}(L)$ & $x_{2}$ & $y_{2}$ & $L^{2}$ \\
\cline{2-3}
$U$ & $x_{1}$ & $y_{1}$ & $L^{3}=Z_{3}(L)$ \\
\cline{2-3}
\end{tabular}  
\mbox{}\\ \\ \\
Notice that $UZ_{4}(L)={\mathbb F}x_{1}y_{2}+{\mathbb F}x_{1}y_{3}=
{\mathbb F}x_{5}+{\mathbb F}x_{1}y_{3}$, where $x_{1}y_{3}\in L^{5}$. 
Again we consider two possible cases. \\ \\
{\bf I. Algebras where $UZ_{4}(L)$ is $1$-dimensional} \\ \\
Here $UZ_{4}(L)=x_{1}Z_{4}(L)={\mathbb F}x_{5}$ and there is a characteristic subspace
$V$ of codimension $1$ in $Z_{4}(L)$ that contains $L^{3}$ given by the formula
    $$V=\{x\in Z_{4}(L):\, Ux=0\}.$$
We can then choose our standard basis such that 
   $$V=L^{3}+{\mathbb F}y_{3}=U+{\mathbb F}y_{1}+{\mathbb F}y_{3}.$$
Notice that in particular $x_{1}y_{3}=0$. From this we also get a 
$1$-dimensional characteristic subspace $V^{2}={\mathbb F}y_{1}y_{3}$. Notice
that $(y_{1}y_{3},y_{2})\not =0$ as otherwise $y_{1}y_{2}\in L^{5}\cap
L^{3}L^{2}$ that contradicts our assumption that $L^{3}L^{2}\not = L^{5}$. Thus
$y_{1}y_{3}\in L^{4}\setminus Z(L)$ and we can choose our standard basis such
that ${\mathbb F}y_{1}y_{3}={\mathbb F}x_{2}$. In fact it is not difficult to see that with the data we have acquired so far we can choose our standard
basis such that 
\begin{equation}
   x_{1}y_{2}=x_{5},\ y_{1}y_{2}=x_{3},\ x_{1}y_{3}=0,\ y_{1}y_{3}=-x_{2}.
\end{equation}
This deals with all triple values apart from 
$$\begin{array}{llll}
   (x_{1}y_{4},y_{5})=a, & (y_{2}y_{3},y_{4})=c, & (y_{2}y_{4},y_{5})=e, &
(x_{2}y_{4},y_{5})=r, \\
  (y_{1}y_{4},y_{5})=b, & (y_{2}y_{3},y_{5})=d, & (y_{3}y_{4},y_{5})=f, &
\mbox{} 
\end{array}$$
Notice that $r\not =0$ as $x_{2}y_{3}=0$ but $x_{2}\not\in Z(L)$. 
We will show that we can choose a new standard basis so that the values of
$a=b=c=d=e=f=0$ and $r=1$. 
Replace $x_1,x_4, y_1, y_2, y_3, y_4, y_5$ by $x_1-(a/r)x_2+c x_4 + d x_5, rx_4, y_1-(b/r)x_2,
y_2-(b/r)x_1 - (e/r)x_2 - (f/r)x_3 - (bc/r)x_4 + (a/r)y_1, y_3-(f/r)x_2 , (1/r)y_4 - (c/r)y_1, y_5 - d y_1+(bd/r)x_2$. Thus we have that we get a unique algebra.
\begin{prop} There is a unique nilpotent SAA $L$ with an isotropic centre of dimension $3$
and where $\mbox{dim\,}L^{3}=6$ that has the further properties that 
$L^{3}L^{2}\not =L^{5}$ and $\mbox{dim\,}UZ_{4}(L)=1$. This algebra is given
by the presentation 
$${\mathcal P}_{10}^{(3,4)}:\ (x_{1}y_{2},y_{5})=1,\ (y_{1}y_{2},y_{3})=1,\ (x_{2}y_{4},y_{5})=1.$$
\end{prop}
{\bf Remark}. As before, inspection shows that the algebra with the presentation above satisfies all the properties listed. \\ \\ 
{\bf II. Algebras where $UZ_{4}(L)$ is $2$-dimensional} \\ \\
Here we can pick our standard basis such that $UZ_{4}(L)={\mathbb F}x_{5}+
{\mathbb F}x_{4}$. As $L^{3}L^{2}\not =L^{5}$ we know that $(y_{1}y_{2},y_{3}) \neq 0$
and from this one sees that $L^{3}Z_{4}(L)=L^{4}$. Furthermore it is not difficult to see that we can choose our standard basis such that
\begin{equation}
     x_{1}y_{2}=x_{5},\ x_{1}y_{3}=x_{4},\ y_{1}y_{2}=x_{3},\ y_{1}y_{3}=-x_{2}.
\end{equation}
In order to clarify the structure further we are only left with the triple values
$$\begin{array}{llll}
   (x_{1}y_{4},y_{5})=a, & (y_{2}y_{3},y_{4})=c, & (y_{2}y_{4},y_{5})=e, &
(x_{2}y_{4},y_{5})=r, \\
  (y_{1}y_{4},y_{5})=b, & (y_{2}y_{3},y_{5})=d, & (y_{3}y_{4},y_{5})=f. &
\mbox{} 
\end{array}$$
Let $\alpha = b -cr$ and replace
$x_1, x_2, x_3, x_4, x_5,  y_1, y_2, y_3, y_4, y_5 $ by 
$x_{1}-(a/r) x_{2}+(\alpha/r) x_{4} + d x_5$,
$(1/r)x_2$,
$rx_3$,
$(1/r)x_4$,
$rx_5$,
$y_{1}- (b/r) x_{2}$,
$ry_{2} - b x_{1} - e x_2 -f x_3 + a y_1$,
$(1/r)y_{3}-(f/r^2)x_{2}$,
$ry_{4} - \alpha(b/r) x_{2} - \alpha y_{1}$,
$(1/r)y_{5} + d(b/r^2) x_{2} - (d/r) y_{1}$.
One checks readily that we get a new standard basis such that $a=b=c=d=e=f=0$ and $r=1$. Thus we see again that we have a unique algebra.
%
%%%%%%%%%%%%%%%%%%%%
\iffalse
{\color{blue}
By replacing $y_3, y_1, y_5, x_1, y_4, y_2$ by 
$ y_{3}-(f/r)x_{2},
y_{1}- (b/r) x_{2},
y_{5} + d(b/r) x_{2} - d y_{1} ,
x_{1}-(a/r) x_{2}+\alpha x_{4} + d x_5 ,
y_{4} - \alpha(b/r) x_{2} - \alpha y_{1},
y_{2}-(b/r) x_{1} - (e/r) x_2 - (f/r)x_3 + (a/r)y_1$,
%
where $\alpha=c- b/r$. One checks readily that we get a new standard basis such that a=b=c=d=e=f=0.
%
Finally changing $x_{2}$,$x_{3}$,$x_{4}$,$x_{5}$,$y_{2}$,$y_{3}$,$y_{4}$,$y_{5}$
by $(1/r)x_{2}$,$rx_{3}$,$(1/r)x_{4}$,$rx_{5}$,$ry_{2}$,$(1/r)y_{3}$,$ry_{4}$,
$(1/r)y_{5}$
gives us $(x_{2}y_{4},y_{5})=1$ and thus we see again that we have a unique
algebra.}
%
\fi
%%%%%%%%%%%%%

%
\begin{prop} There is a unique nilpotent SAA $L$ with an isotropic centre of dimension $3$
and where $\mbox{dim\,}L^{3}=6$ that has the further properties that 
$L^{3}L^{2}\not =L^{5}$ and $\mbox{dim\,}UZ_{4}(L)=2$. This algebra is given
by the presentation 
$${\mathcal P}_{10}^{(3,5)}:\ (x_{1}y_{2},y_{5})=1,\ (y_{1}y_{2},y_{3})=1,\ 
(x_{1}y_{3},y_{4})=1,\ (x_{2}y_{4},y_{5})=1.$$
\end{prop}
\subsubsection{The algebras where $L^{3}L^{2}=L^{5}$}
$$\begin{array}{ll}
\begin{tabular}{c|c|c|c} 
\cline{2-3} 
 \mbox{} & $x_{5}$ & $y_{5}$ & \mbox{} \\
   $L^{5}$ & $x_{4}$ & $x_{4}$ & \mbox{} \\
\cline{2-3}
\mbox{}$Z(L)$  & $x_{3}$ & $y_{3}$ & $Z_{4}(L)$ \\ 
\cline{2-3} 
 \mbox{}$L^{4}=Z_{2}(L)$ & $x_{2}$ & $y_{2}$ & $L^{2}$ \\
\cline{2-3}
\mbox{} & $x_{1}$ & $y_{1}$ & $L^{3}=Z_{3}(L)$ \\
\cline{2-3}
\end{tabular} &  
\begin{array}{l}
   L^{5}={\mathbb F}x_{5}+{\mathbb F}x_{4} \\
 Z(L)={\mathbb F}x_{5}+{\mathbb F}x_{4}+{\mathbb F}x_{3} \\ 
 Z_{2}(L)=L^{4}={\mathbb F}x_{5}+{\mathbb F}x_{4}+{\mathbb F}x_{3}+
{\mathbb F}x_{2} \\
 Z_{3}(L)=L^{3}=Z(L)+{\mathbb F}x_{2}+{\mathbb F}x_{1}+{\mathbb F}y_{1} \\
 L^{2}=L^{3}+{\mathbb F}y_{2} \\
  Z_{4}(L)=L^{3}+{\mathbb F}y_{2}+{\mathbb F}y_{3}
\end{array}
\end{array}$$
Here we are assuming that $L^{5}=L^{3}L^{2}={\mathbb F}x_{1}y_{2}+
{\mathbb F}y_{1}y_{2}$ and thus in particular we know that 
$x_{1}y_{2},y_{1}y_{2}$ is a basis for $L^{5}$. We will now introduce
some linear maps that will help us in understanding the structure. Consider
first the linear maps
$$\begin{array}{ll}
  \phi:\,L^{3}/L^{4}\rightarrow L^{5},\ \bar{u}=u+L^{4}\mapsto u\cdot y_{2} \\
  \psi:\,L^{3}/L^{4}\rightarrow L^{5},\ \bar{u}=u+L^{4}\mapsto u\cdot y_{3}.
\end{array}$$
As $L^{4}Z_{4}(L)=\{0\}$, these maps are well defined. As $L^{3}L^{2}=L^{5}$ 
we also know that $\phi$ is bijective. We thus have the linear map
   $$\tau=\psi\phi^{-1}:\,L^{5}\rightarrow L^{5}.$$
It is the map $\tau$ that will be our key towards understanding the
structure of the algebra. 
\begin{lemm} The minimal polynomial of $\tau=\psi\phi^{-1}$ must be of
degree $2$.
\end{lemm}
{\bf Proof\ \ }We argue by contradiction and
suppose that $\tau=\lambda \mbox{id}$. Replacing $y_{3},x_{2}$ by
$y_{3}-\lambda y_{2}, x_{2}+\lambda x_{3}$ gives us a new standard basis
where $\tau=0$. Pick our standard basis such that $\bar{x}_{1}=x_{1}+
L^{4}=\phi^{-1}(x_{4})$ and $\bar{y}_{1}=y_{1}+L^{4}=\phi^{-1}(x_{5})$. We then
have
  $$x_{1}y_{2}=x_{4},\ y_{1}y_{2}=x_{5},\ x_{1}y_{3}=0,\ y_{1}y_{3}=0.$$
Now $y_{2}y_{3}\perp x_{1},y_{1},y_{2},y_{3}$ and thus 
      $$y_{2}y_{3}=a x_{4}+b x_{5}$$
for some $a,b\in {\mathbb F}$. Replacing $y_{3}$ by $y_{3}+ax_{1}+by_{1}$,
$x_{1}$ by $x_{1}-bx_{3}$ and $y_{1}$ by $y_{1}+ax_{3}$, we can assume that
$y_{2}y_{3}=0$. \\ \\
Now suppose that $(y_{3}y_{4},y_{5})=a$ and $(x_{2}y_{4},y_{5})=b$. Notice
that $b\not =0$ as $x_{2} \not \in Z(L)$ and $x_{2}y_{3}=0$. Replace $y_{3},y_{2}$
by $y_{3}-(a/b)x_{2},\ y_{2}-(a/b)x_{3}$ and we get a new standard basis
where all the previous identities hold but also $(y_{3}y_{4},y_{5})=0$. We
thus get the contradiction that $y_{3}\in Z(L)$. $\Box$\\ \\
Notice next that if we have an alternative
standard basis $\tilde{x}_{1},\tilde{x}_{2},
\ldots ,\tilde{y}_{5}$, then $\tilde{y}_{2}=cy_{2}+u$ and $\tilde{y}_{3}=
ay_{3}+by_{2}+v$ where $a, c\not =0$ and where $u,v\in L^{3}$. If the minimal polynomial of $\tau$ 
with respect to the old basis is $f(t)$ then the minimal polynomial with
respect to the new basis is a multiple of $f((c/a)(t-(b/c))$. In particular
we have the following possible distinct scenarios that do not depend on what 
standard basis we choose. \\ \\
{\bf A}. The minimal polynomial of $\tau$ has two distinct roots in ${\mathbb F}$. \\ 
{\bf B}. The minimal polynomial of $\tau$ has a double root in ${\mathbb F}$  \\ 
{\bf C}. The minimal polynomial of $\tau$ is irreducible in ${\mathbb F}[t]$. \\ \\
\large
{\bf I. Algebras of type A}. \\ \\
\normalsize
Suppose the two distinct roots of the minimal polynomial of $\tau=\psi\phi^{-1}$ are $\lambda$ and
$\mu$. Pick some eigenvectors $x_{4}$ and $x_{5}$ with respect to the
eigenvalues $\lambda$ and $\mu$ respectively. Thus 
\begin{eqnarray*}
       \psi\phi^{-1}(x_{4}) & = & \lambda x_{4}, \\
       \psi\phi^{-1}(x_{5}) & = & \mu x_{5}.
\end{eqnarray*}
Replacing $y_{3}, x_{2}$ by $y_{3}-\lambda y_{2}, x_{2}+\lambda x_{3}$ we see
that $\psi\phi^{-1}(x_{4})=0$ and we can thus assume that $\lambda =0$. Then
replace $y_{3},x_{3}$ by $(1/\mu)y_{3},\mu x_{3}$ and we get that 
$\psi\phi^{-1}(x_{5})=x_{5}$ and we can now assume that $\mu =1$. \\ \\ 
We would like to pick our standard basis such that $\bar{x}_{1}=x_{1}+L^{4}=
\phi^{-1}(x_{4})$ and $\bar{y}_{1}=y_{1}+L^{4}=
\phi^{-1}(x_{5})$. The only problem here is that we
need $(x_{1},y_{1})=1$ but this
can be easily arranged. If $(x_{1},y_{1})=\sigma$ then we just need to replace
$y_{1},x_{5},y_{5}$ by $(1/\sigma)y_{1},(1/\sigma)x_{5},\sigma y_{5}$. We have thus seen that we can choose our standard basis such that 
\begin{equation}
 x_{1}y_{2}=x_{4},\ y_{1}y_{2}=x_{5},\ x_{1}y_{3}=0,\ y_{1}y_{3}=x_{5}.
\end{equation}
Recall also that $x_{2}y_{3}=0$ since $L^{4}Z_{4}(L)=\{0\}$. In order to fully
determine the structure of the algebra we are only left with the following
triple values
$$\begin{array}{llll}
  (x_{1}y_{4},y_{5})=a, & (y_{2}y_{3},y_{4})=c, & (y_{2}y_{4},y_{5})=e, &
                    (x_{2}y_{4},y_{5})=r, \\
 (y_{1}y_{4},y_{5})=b, & (y_{2}y_{3},y_{5})=d, & (y_{3}y_{4},y_{5})=f. &
\end{array}$$
Notice that $r\not =0$ as $x_{2}y_{3}=0$ but $x_{2}\not\in Z(L)$. 
Let $\alpha=(1/r)(e - f - bd  + ab/r - ac)$ and
replace $y_4, y_3, y_2, y_1, x_1$ by
$y_4 + e y_1 - e(b/r) x_2  + ce x_3,
y_3 + c x_1 - \alpha x_2 + (d+a/r)y_1,
y_2 + (a/r) y_1 - (b/r) x_1 + ((e-f)/r) x_3,
y_1 - (b/r)x_2+c x_3,
x_1 -(a/r)x_2  - (d+a/r) x_3 - e x_4$
One checks readily that we get a new standard basis such that $a=b=c=d=e=f=0$.
We have thus seen that $L$ has a presentation of the form ${\mathcal P}_{10}^{(3,6)}(r)$ as described in the next proposition.
\begin{prop} Let $L$ be a nilpotent SAA of dimension $10$ with an isotropic centre of dimension $3$ that has the further properties that $\mbox{dim\,}L^{3}=6$, $L^{3}L^{2}=L^{5}$ and $L$ is of type A. Then $L$ has a presentation of the form 
   $${\mathcal P}_{10}^{(3,6)}(r):\ (x_{2}y_{4},y_{5})=r,\ (x_{1}y_{2},y_{4})=1,\ (y_{1}y_{2},y_{5})=1,\
	    (y_{1}y_{3},y_{5})=1$$
where $r\not =0$. Furthermore the presentations ${\mathcal P}_{10}^{(3,6)}(r)$ and ${\mathcal P}_{10}^{(3,6)}(s)$
describe the same algebra if and only if $s/r\in ({\mathbb F}^{*})^{3}$. 
\end{prop}
{\bf Proof}\ \ We have already seen that all such algebras have a presentation of the form ${\mathcal P}_{10}^{(3,6)}(r)$ for some $0\not =r\in {\mathbb F}$. Straightforward calculations show that conversely any algebra with such a presentation has the properties stated in the Proposition. It remains to prove the isomorphism property. To see that the property is sufficient, suppose we have an algebra $L$ with presentation ${\mathcal P}_{10}^{(3,6)}(r)$ with respect to some given standard basis. Let $s$ be any element in ${\mathbb F}^{*}$ such that $s/r=b^{3}\in 
({\mathbb F}^{*})^{3}$. Replace the basis with a new standard basis $\tilde{x}_{1},\ldots ,\tilde{y}_{5}$ where $\tilde{x}_{1}=x_{1}$, $\tilde{y}_{1}=y_{1}$, $\tilde{x}_{2}=bx_{2}$, 
$\tilde{y}_{2}=(1/b)y_{2}$, $\tilde{x}_{3}=bx_{3}$, 
$\tilde{y}_{3}=(1/b)y_{3}$, $\tilde{x}_{4}=(1/b)x_{4}$, 
$\tilde{y}_{4}=by_{4}$, $\tilde{x}_{5}=(1/b)x_{5}$, 
$\tilde{y}_{5}=by_{5}$. Direct calculations show that $L$ has 
presentation ${\mathcal P}_{10}^{(3,6)}(s)$ with respect to the new basis. \\ \\
It remains to see that the property is necessary. Consider again an algebra $L$ with presentation ${\mathcal P}_{10}^{(3,6)}(r)$ and suppose that $L$ has also a presentation ${\mathcal P}_{10}^{(3,6)}(s)$ with respect to some other standard basis
$\tilde{x}_{1},\ldots ,\tilde{y}_{5}$. We want to show that $s/r\in ({\mathbb F}^{*})^{3}$. We know that $L={\mathbb F}\tilde{y}_{5}+{\mathbb F}\tilde{y}_{4}+Z_{4}(L)={\mathbb F}y_{5}+{\mathbb
F}y_{4}+Z_{4}(L)$. Thus 
\begin{eqnarray*}
       \tilde{y}_{4} & = & ay_{4}+by_{5}+u \\
    \tilde{y}_{5} & = & cy_{4}+dy_{5}+v
\end{eqnarray*}
for some $u,v\in Z_{4}(L)$ and $a,b,c,d\in {\mathbb F}$ where
$ad-bc\not =0$. As $L^{3}L^{2}=L^{5}\perp Z_{4}(L)$ and as
$Z_{4}(L)L^{4}=0$ we have $(Z_{4}(L)L^{2},L^{3})=(Z_{4}(L)L,L^{4})=0$ and thus $Z_{4}(L)L^{2}\leq (L^{3})^{\perp}=L^{4}$ and $Z_{4}(L)L\leq (L^{4})^{\perp}=L^{3}$. It follows that 
\begin{eqnarray*}
   \tilde{y}_{4}\tilde{y}_{5}\tilde{y}_{5} & = & 
	  (ay_{4}+by_{5})(cy_{4}+dy_{5})(cy_{4}+dy_{5})+w \\
   \tilde{y}_{5}\tilde{y}_{4}\tilde{y}_{4} & = & 
	(cy_{4}+dy_{5})(ay_{4}+by_{5})(ay_{4}+by_{5}) +z
\end{eqnarray*}
for some $w,z\in L^{4}$. Using the fact that $(L^{4},L^{3})=0$,
as $L^{6}=0$, we then see that 
      $$s^{2}=(\tilde{y}_{4}\tilde{y}_{5}\tilde{y}_{5},
			   \tilde{y}_{5}\tilde{y}_{4}\tilde{y}_{4})=
				r^{2}(ad-bc)^{3}.$$
Hence $s/r\in ({\mathbb F}^{*})^{3}$. $\Box$ \\ \\
{\bf Remark}. Notice that it follows that we have only one 
algebra if $({\mathbb F}^{*})^{3}={\mathbb F}^{*}$. This includes all fields that are algebraically closed as well as 
${\mathbb R}$. For a finite field of order $p^{n}$ there are
$3$ algebras if $3|p^{n}-1$ but otherwise one. For ${\mathbb Q}$ there are infinitely many algebras. \\ \\ 
\large
{\bf II. Algebras of type B}. \\ \\
\normalsize
Suppose that the double root of the minimal polynomial of $\tau=\psi\phi^{-1}$ is $\lambda$. We can then have  a
basis $x_{4},x_{5}$ for $L^{5}$ such that 
\begin{eqnarray*}
    \psi\phi^{-1}(x_{4}) & = & \lambda x_{4} \\
   \psi\phi^{-1}(x_{5}) & = & \lambda x_{5}+x_{4}.
\end{eqnarray*}
If we replace $y_{3},x_{2}$ by $y_{3}-\lambda y_{2}, x_{2}+\lambda x_{3}$ then
we can furthermore assume that $\lambda =0$. We want to pick our standard
basis such that $\bar{x}_{1}=x_{1}+L^{4}=\phi^{-1}(x_{4})$ and $\bar{y}_{1}=
y_{1}+L^{4}$. Again the only problem is to arrange for $(x_{1},y_{1})=1$. But
if $(x_{1},y_{1})=\sigma$ then we replace $x_{5},x_{3},y_{1},y_{3},y_{5}$ by
$(1/\sigma)x_{5},(1/\sigma)x_{3}, (1/\sigma)y_{1},\sigma y_{3},\sigma y_{5}$ and that gives $(x_{1},y_{1})=1$. We have thus seen that we can choose our standard
basis such that 
\begin{equation}
    x_{1}y_{2}=x_{4},\ y_{1}y_{2}=x_{5},\ x_{1}y_{3}=0,\ y_{1}y_{3}=x_{4}.
\end{equation}
As before we have furthermore $x_{2}y_{3}=0$ and we are only left with the
following triple values 
$$\begin{array}{llll}
  (x_{1}y_{4},y_{5})=a, & (y_{2}y_{3},y_{4})=c, & (y_{2}y_{4},y_{5})=e, &
                              (x_{2}y_{4},y_{5})=r, \\
  (y_{1}y_{4},y_{5})=b, & (y_{2}y_{3},y_{5})=d, & (y_{3}y_{4},y_{5})=f.
\end{array}$$
Notice that $r\not =0$ as $x_{2}y_{3}=0$ but $x_{2}\not\in Z(L)$. 
Let $ \alpha =( a^2+acr+bdr )/r$, $\beta = cr+a$ and replace
 $x_1, y_3, y_1, y_4, y_2, y_5$ by
$x_1 -(a/r) x_2 - d x_3 - e x_4 + f x_5 ,
y_3 - $
$(\alpha/r) x_2  + (\beta/r) x_1+ d y_1  ,
y_1 - (b/r) x_2 + (\beta/r) x_3 ,
y_4 - e(b/r) x_2  + e(\beta/r) x_3 + e y_1  ,
y_2 + (a/r)y_1 - (b/r) x_1,
y_5 + f (b/r) x_2 - f (\beta/r) x_3  - f y_1  $
One checks readily that these changes do not affect $(17)$ and we get a new standard basis where $a=b=c=d=e=f=0.$
We thus arrive at at presentation of the form ${\mathcal P}_{10}^{(3,7)}(r)$ as given
in the next proposition. 
\begin{prop} Let $L$ be a nilpotent SAA of dimension $10$ with an isotropic centre of dimension $3$ that has the further properties that $\mbox{dim\,}L^{3}=6$, $L^{3}L^{2}=L^{5}$ and $L$ is of type B. Then $L$ has a presentation of the form 
   $${\mathcal P}_{10}^{(3,7)}(r):\ (x_{2}y_{4},y_{5})=r,\ (x_{1}y_{2},y_{4})=1,\ (y_{1}y_{2},y_{5})=1,\
	    (y_{1}y_{3},y_{4})=1$$
where $r\not =0$. Furthermore the presentations ${\mathcal P}_{10}^{(3,6)}(r)$ and ${\mathcal P}_{10}^{(3,6)}(s)$
describe the same algebra if and only if $s/r\in ({\mathbb F}^{*})^{3}$. 
\end{prop}
{\bf Proof}\ \ Similar to the proof of Proposition $4.11$. \\ \\
\large
{\bf III. Algebras of type C}. \\ \\
\normalsize
It turns out to be useful to consider the cases ${\mbox char\,}{\mathbb F}\not =2$ and 
${\mbox char\,}{\mathbb F}=2$ separately. \\ \\
\underline{{\bf a}. The algebras where ${\mbox  char\,}{\mathbb F}\not=2$} \\ \\
Suppose the minimal polynomial of $\tau=\psi\phi^{-1}$ is $t^{2}+at+b$ with respect to some $y_{2},y_{3}$. Replacing $y_{3}$ by $y_{3}+(a/2)y_{2}$, one gets a minimal polynomial of the form $t^{2}-s$ with
$s\not\in {\mathbb F}^{2}$. \\ \\
{\bf Remark}. Let $\tilde{y}_{3}=\alpha y_{3}+u$ where $\alpha\not =0$ and 
$u\in L^{2}$. For the minimal polynomial of $\tau$ to have trivial linear term
we must have $u\in L^{3}$. Thus ${\mathbb F}y_{3}+L^{3}$ is a characteristic
subspace of $L$. \\ \\
Pick any $0\not =x_{5}\in L^{5}$ and let $x_{4}=\psi\phi^{-1}(x_{5})$. Then
$\psi\phi^{-1}(x_{4})=sx_{5}$. We want to pick our standard basis such that $\phi^{-1}(x_{4})=x_{1}+L^{4},
\phi^{-1}(x_{5})=y_{1}+L^{4}$. For this to work out we need $(x_{1},y_{1})=1$. Again this can be easily arranged. If $(x_{1},y_{1})=\sigma$, then we replace $x_{5},y_{1},y_{3}$ by $(1/\sigma)x_{5}$, 
$(1/\sigma)y_{1}$, $\sigma y_{3}$ and we get $(x_{1},y_{1})=1$ and $\psi\phi^{-1}(x_{4})=
(\sigma^{2}s)x_{5}$. We have thus seen that we choose our standard basis such that 
\begin{equation}
   x_{1}y_{2}=x_{4},\ y_{1}y_{2}=x_{5},\ x_{1}y_{3}=sx_{5},\ y_{1}y_{3}=x_{4}
\end{equation}
for some $s\not\in {\mathbb F}^{2}$. In order to clarify the structure further we are only left with the following triple values 
$$\begin{array}{llll} 
  (x_{1}y_{4},y_{5})=a, & (y_{2}y_{3},y_{4})=c, & (y_{2}y_{4},y_{5})=e, & (x_{2}y_{4},y_{5})=r, \\
  (y_{1}y_{4},y_{5})=b, & (y_{2}y_{3},y_{5})=d, & (y_{3}y_{4},y_{5})=f.
\end{array}$$
Notice that $r\not =0$ as $x_{2}y_{3}=0$ but $x_{2}\not\in Z(L)$. 
Let $\alpha = (a^2 + acr - b^2 s+ bdr)r $, $\beta = c+ a/r$, $\gamma =(s(b/r)-d ) $ and replace
$x_1, y_4,y_1, y_5, y_2, y_3$ by
$x_1 - (a/r) x_2 - e x_4+ f x_5   - \gamma x_3  ,
y_4 - e(b/r) x_2 + e \beta x_3 + e y_1 ,
y_1 - (b/r) x_2+ \beta x_3   ,
y_5 + f (b/r) x_2 - f \beta x_3  - f y_1,
y_2 - (b/r) x_1 + (a/r) y_1,
y_3 - (\alpha/r) x_2 + \beta x_1 - \gamma y_1 ,
$
We see then that the remaining triple values are zero.
Thus $L$ has a presentation of the form 
${\mathcal P}_{10}^{(3,8)}(r,s)$ as described in the next proposition. 
\begin{prop} Let $L$ be a nilpotent SAA of dimension $10$ over a field of characteristic that is not $2$ that has an isotropic centre of dimension $3$. Suppose also that $L$ has the further properties that $\mbox{dim\,}L^{3}=6$, $L^{3}L^{2}=L^{5}$ and $L$ is of type C. Then $L$ has a presentation of the form 
\begin{align*}
{\mathcal P}_{10}^{(3,8)}(r,s):\ &(x_{2}y_{4},y_{5})=r,\ (x_{1}y_{3},y_{5})=s,\ (x_{1}y_{2},y_{4})=1,\ (y_{1}y_{2},y_{5})=1,\ \\
	   & (y_{1}y_{3},y_{4})=1,
\end{align*}
where $r\not =0$ and $s\not\in {\mathbb F}^{2}$. Furthermore the presentations ${\mathcal P}_{10}^{(3,8)}(\tilde{r},\tilde{s})$ and ${\mathcal P}_{10}^{(3,8)}(r,s)$
describe the same algebra if and only if $\frac{\tilde{r}}{r}\in ({\mathbb F}^{*})^{3}$ and 
$\frac{s}{\tilde{s}}\in G(s)$ where $G(s)=\{(x^{2}-y^{2}s)^{2}:\ (x,y)\in {\mathbb F}\times {\mathbb F}\setminus \{(0,0)\}\}$. 
\end{prop}
{\bf Proof}\ \ We have already seen that any such algebra has a presentation of the given form. Direct calculations show that an algebra with a presentation ${\mathcal P}_{10}^{(3,8)}(r,s)$ has the properties stated. We turn to the isomorphism property. To see that the condition is sufficient, suppose we have an algebra $L$ that has presentation ${\mathcal P}_{10}^{(3,8)}(r,s)$ with respect to some standard basis $x_{1},y_{1},\ldots ,x_{5},y_{5}$. Suppose that $\frac{\tilde{r}}{r}=\frac{1}{\beta^{3}}$
and $\frac{s}{\tilde{s}}=[(b/\beta)^{2}-s(a/\beta)^{2}]^{2}$
for some $\beta \in \mathbb{F} \setminus \{ 0\}$ and $(a,b) \in {\mathbb F}\times {\mathbb F}\setminus
\{(0,0)\}$. Let $\alpha=\frac{\beta^{3}}{b^{2}-a^{2}s}$ and consider a new standard basis
$$\begin{array}{ll}
   \tilde{x}_{1}=(\alpha/\beta^{2})(bx_{1}+asy_{1}), & 
	              \tilde{y}_{1}=(1/\beta)(by_{1}+ax_{1}), \\
  \tilde{x}_{2}=(1/\beta)x_{2}, & \tilde{y}_{2}=\beta y_{2} \\
  \tilde{x}_{3}=(1/\alpha)x_{3}, & \tilde{y}_{3}=\alpha y_{3} \\
  \tilde{x}_{4}=(\alpha/\beta)(bx_{4}+asx_{5}), & \tilde{y}_{4}=(1/\beta^{2})(by_{4}-ay_{5}), \\
 \tilde{x}_{5}=ax_{4}+bx_{5}, & \tilde{y}_{5}=(\alpha/\beta^{3})(by_{5}-asy_{4}),
\end{array}$$
Calculations show that $L$ has then presentation 
${\mathcal P}_{10}^{(3,8)}(\tilde{r},\tilde{s})$ with respect to the new standard basis. \\ \\
It remains to see that the conditions are also necessary. Consider an algebra 
$L$ with presentation ${\mathcal P}_{10}^{(3,8)}(r,s)$ with respect to some standard basis $x_{1},y_{1},\ldots ,x_{5},y_{5}$. Take some arbitrary new standard 
basis $\tilde{x}_{1},\tilde{y}_{1},\ldots ,\tilde{x}_{5},\tilde{y}_{5}$ such
that $L$ satisfies the presentation ${\mathcal P}_{10}^{(3,8)}(\tilde{r},
\tilde{s})$ with respect to the new basis for some $0\neq \tilde{r}\in 
{\mathbb F}$ and $\tilde{s} \not \in {\mathbb F}^{2}$. Notice that 
\begin{eqnarray*}
        \tilde{x}_{5} & = & ax_{4}+bx_{5} \\
        \tilde{y}_{2} & = & \beta y_{2} + u \\
        \tilde{y}_{3} & = & \gamma y_{3} + v,
\end{eqnarray*}
such that $u,v\in L^{3}$ and $0\not =\alpha,\beta \in {\mathbb F}$. The reader
can convince himself that $\tilde{r}/r\in ({\mathbb F}^{*})^{3}$ and 
$s/\tilde{s}\in G(s)$. $\Box$ \\ \\
{\bf Examples}. (1) If ${\mathbb F}={\mathbb C}$ then as any quadratic polynomial is reducible, there are not algebras of type C. This holds more generally for any field ${\mathbb F}$ of characteristic that is not $2$ and where all the elements in
${\mathbb F}$ have a square root in ${\mathbb F}$. \\ \\
(2) Suppose ${\mathbb F}={\mathbb R}$. Let $s\not \in {\mathbb R}^{2}$ and
$0\not =r\in {\mathbb R}$. Then $1/r\in ({\mathbb R}^{*})^{3}$ and 
$s<0$. Also $s/(-1)=a^{4}=(a^{2}-0^{2}s)^{2}$ for some $a\in {\mathbb R} \setminus \{0\}$.
We thus have that ${\mathcal P}_{10}^{(3,8)}(r,s)$ describes the same algebra
as ${\mathcal P}_{10}^{(3,8)}(1,-1)$. There is thus a unique algebra in this
case. \\ \\
(3) Let ${\mathbb F}$ be a finite field of some odd characteristic $p$. Suppose
that $|{\mathbb F}|=p^{n}$. Let $s$ be any element that is not in 
$({\mathbb F}^{*})^{2}$. Notice then that ${\mathbb F}^{*}=({\mathbb F}^{*})^{2}\cup s({\mathbb F}^{*})^{2}$ and thus for any $\tilde{s}$ that is not
in ${\mathbb F}^{2}$, we have $s/\tilde{s}\in ({\mathbb F}^{*})^{2}=G(s)$. 
We can thus keep $s$ fixed and each algebra has a presentation of
the form ${\mathcal Q}(r)={\mathcal P}_{10}^{(3,8)}(r,s)$ where ${\mathcal Q}(\tilde{r})$ and ${\mathcal Q}(r)$ describe the same algebra if and only
if $\tilde{r}/r\in ({\mathbb F}^{*})^{3}$. There are thus either three or
one algebra according to whether $3$ divides $p^{n}-1$ or not. $\Box$ \\ \\
\underline{{\bf b}. The algebras where ${\mbox  char\,}{\mathbb F}=2$} \\ \\
If the irreducible minimal polynomial of $\psi\phi^{-1}$ is $t^{2}+rt+s$
with respect to $y_{2},y_{3}$ then the minimal polynomial with respect to
$ay_{2},by_{3}+cy_{2}$, where $a,b\not =0$, is 
   $$t^{2}+r(b/a)t+[(c/a)^{2}+r(c/a)(b/a)+(b/a)^{2}s].$$
Thus we have two distinct subcases (that do not depend on the choice of 
basis). Let $m=m(y_{2},y_{3})$ be the minimal polynomial of $\psi\phi^{-1}$
with respect to a given standard basis for $L$.  \\ \\
(1) The minimal polynomial $m$  is of the form $t^{2}-s$ for
some $s\not\in{\mathbb F}^{2}$. \\ \\
(2) The minimal polynomial $m$ is of the form $t^{2}+rt+s$ where $r\not =0$ and the polynomial is irreducible. \\ \\
For case (1) we get the same situation as in Proposition 4.13. 
\begin{prop} Let $L$ be a nilpotent SAA of dimension $10$ over a field of characteristic $2$ that has an isotropic centre of dimension $3$. Suppose also that $L$ has the further properties that $\mbox{dim\,}L^{3}=6$, $L^{3}L^{2}=L^{5}$ and $L$ is of type C where the minimal polynomial $m(y_{2},y_{3})$ is of the
form $t^{2}-s$ for some $s\not\in {\mathbb F}^{2}$.  Then $L$ has a presentation of the form 
\begin{align*}
    {\mathcal P}_{10}^{(3,8)}(r,s):\ &(x_{2}y_{4},y_{5})=r,\ (x_{1}y_{3},y_{5})=s,\ (x_{1}y_{2},y_{4})=1,\ (y_{1}y_{2},y_{5})=1,\ \\
	    &(y_{1}y_{3},y_{4})=1,
\end{align*}
where $r\not =0$ and $s\not\in {\mathbb F}^{2}$. Furthermore the presentations ${\mathcal P}_{10}^{(3,8)}(\tilde{r},\tilde{s})$ and ${\mathcal P}_{10}^{(3,8)}(r,s)$
describe the same algebra if and only if $\frac{\tilde{r}}{r}\in ({\mathbb F}^{*})^{3}$ and 
$\frac{s}{\tilde{s}}\in G(s)$ where $G(s)=\{(x^{2}-y^{2}s)^{2}:\ (x,y)\in {\mathbb F}\times {\mathbb F}\setminus \{(0,0)\}\}$. 
\end{prop}
{\bf Example}. Consider the field ${\mathbb Z}_{2}(x)$ of rational functions in
one variable over ${\mathbb Z}_{2}$. Notice that 
$${\mathbb Z}_{2}(x)^{*}=\{f(x)^{2}+xg(x)^{2}:\,f(x),g(x)\in {\mathbb Z}_{2}(x)^{2}\setminus \{(0,0)\}\}.$$
Thus $G(x)=({\mathbb Z}_{2}(x)^{*})^{2}$ and the last proposition tells us that 
${\mathcal P}_{10}^{(3,8)}(\tilde{r}(x),\tilde{s}(x))$ and 
${\mathcal P}_{10}^{(3,8)}(r(x),s(x))$ describe the same algebra if and only
if $\tilde{r}(x)/r(x)\in ({\mathbb Z}_{2}(x)^{*})^{3}$ and 
$s(x)/\tilde{s}(x)\in ({\mathbb Z}_{2}(x)^{*})^{2}$. We thus have infinitely
many algebras here. \\ \\
We then move on to the latter collection of algebras. For the rest of 
this subsection we will be assuming that the minimal polynomial of $\psi\phi^{-1}$ is an irreducible polynomial of the form $t^{2}+rt+s$ where $r\not =0$. \\ \\
Pick $0\not =x_{5}\in L^{5}$ and let $x_{4}=\psi\phi^{-1}(x_{5})$. Then $\psi\phi^{-1}(x_{4})=rx_{4}+sx_{5}$. We want to pick our standard basis such that 
$x_{1}+L^{4}=\phi^{-1}(x_{4})$ and $y_{1}+L^{4}=
\phi^{-1}(x_{5})$. The only constraint to worry about is, as before, that
$(x_{1},y_{1})=1$. If $(x_{1},y_{1})=\sigma$, we just need to then replace
$y_{3}$ by $(1/\sigma) y_{3}$. Notice that this changes the minimal
polynomial of $\psi\phi^{-1}$ to $t^{2}+(r/\sigma)t+(s/\sigma^{2})$. In any case 
this shows that we can choose our standard basis such that 
\begin{equation}
      x_{1}y_{2}=x_{4},\ y_{1}y_{2}=x_{5},\ x_{1}y_{3}=rx_{4}+sx_{5},\ y_{1}y_{3}=x_{4}
\end{equation}
for some $r,s\in {\mathbb F}$ where $r\not =0$ and $t^{2}+rt+s$ is irreducible. 
As before we also know that $x_{2}y_{3}=0$. In order to clarify the structure 
further we are only left with the following triple values:
$$\begin{array}{llll}
   (x_{1}y_{4},y_{5})=a, & (y_{2}y_{3},y_{4})=c, & (y_{2}y_{4},y_{5})=e, & (x_{2}y_{4},y_{5})=\alpha, \\
   (y_{1}y_{4},y_{5})=b, & (y_{2}y_{3},y_{5})=d, & (y_{3}y_{4},y_{5})=f. & 
\end{array}$$
Notice that $\gamma \not =0$ as $x_{2}y_{3}=0$ but $x_{2}\not\in Z(L)$. 
Let $z_1=  (a+ c \gamma-br)/\gamma$, $z_2 = (d - s(b/\gamma)) $ and
$z_3=  z_1(a/\gamma)  + z_2 (b/\gamma)$
and replace $x_1, y_3, y_1, y_4, y_2, y_5$ by
$x_1 - (a/\alpha) x_2  - z_2 x_3 - e x_4 + f x_5,
y_3 + z_2 y_1 + z_1   x_1   - z_3 x_2,
y_1 - (b/\gamma) x_2 + z_1 x_3,
y_4 - e(b/\gamma)  x_2 + e y_1   + e z_1 x_3,
y_2 + (a/\gamma) y_1 - (b/\gamma) x_1,
y_5  + f(b/\gamma) x_2 - f  y_1 - f z_1 x_3$.
This
show that we can choose a new standard basis so that the remaining values are
zero. 
We have thus arrived at a presentation of the form ${\mathcal P}_{10}^{(3,9)}$
as described in next proposition. Before stating that proposition we introduce
two groups that are going to play a role. \\ \\
{\bf Definition}. For each minimal polynomial $t^{2}+rt+s$, we let
\begin{eqnarray*}
     H(r) & = & \{x^{2}+rx:\,x\in {\mathbb F}\} \\
   G(r,s) & = & \{x^{2}+rxy+sy^{2}:\, (x,y)\in {\mathbb F}\times 
    {\mathbb F}\setminus \{(0,0)\}\}.
\end{eqnarray*}
{\bf Remarks}. (1) $H(r)$ is a subgroup of the additive group of ${\mathbb F}$. 
\\ \\
(2) Consider the splitting field ${\mathbb F}[\alpha]$ of the polynomial 
$t^{2}+rt+s$ in ${\mathbb F}[t]$. Then $a^{2}+abr+b^{2}s$ is the norm
$N(a+b\alpha)=(a+b\alpha)(a+b(\alpha +r))$ of $a+b\alpha$. As this is a multiplicative function we have that $G(r,s)$ is a multiplicative subgroup of
${\mathbb F}^{*}$. 
\begin{prop} Let $L$ be a nilpotent SAA of dimension $10$ over a field of characteristic $2$ that has an isotropic centre of dimension $3$. Suppose also that $L$ has the further properties that $\mbox{dim\,}L^{3}=6$, $L^{3}L^{2}=L^{5}$ and $L$ is of type C where the minimal polynomial $m(y_{2},y_{3})$ is irreducible with a non-zero linear term. Then $L$ has a presentation of the form 
\begin{eqnarray*}
   {\mathcal P}_{10}^{(3,9)}(\gamma, r,s): &  &  (x_{2}y_{4},y_{5})=\gamma,\ 
(x_{1}y_{3},y_{4})=r,\ (x_{1}y_{3},y_{5})=s,\ \\
            & & (x_{1}y_{2},y_{4})=1,\ (y_{1}y_{2},y_{5})=1,\, (y_{1}y_{3},y_{4})=1
\end{eqnarray*}
where $\gamma, r \not =0$ and $t^{2}+rt+s$ is irreducible. Furthermore the presentations ${\mathcal P}_{10}^{(3,9)}(\tilde{\gamma}, \tilde{r},\tilde{s})$ and ${\mathcal P}_{10}^{(3,9)}(\gamma, r,s)$
describe the same algebra if and only if $\frac{\tilde{\gamma}}{\gamma}\in ({\mathbb F}^{*})^{3}$, $\frac{\tilde{r}}{r}\in G(r,s)$ and $\tilde{s}-(\frac{\tilde{r}}{r})^{2}s\in H(\tilde{r})$. 
\end{prop}
{\bf Proof}\ \ We have already seen that any such algebra has a presentation of the given form. Direct calculations show that conversely any algebra with
a presentation of this type satisfies all the properties listed. It remains to
deal with the isomorphism property. To see that the condition is sufficient, suppose we have an algebra $L$ that has a presentation ${\mathcal P}_{10}^{(3,9)}(\gamma, r,s)$ with respect to some standard basis 
$x_{1},y_{1},\ldots ,x_{5},y_{5}$. Suppose that $\frac{\tilde{\gamma}}{\gamma}=
\frac{1}{\beta^{3}}$, $\frac{r}{\tilde{r}}=(\frac{b}{\beta})^{2}+
(\frac{b}{\beta})(\frac{a}{\beta})r+(\frac{a}{\beta})^{2}s$
and $\tilde{s}-(\frac{\tilde{r}}{r})^{2}s=(\frac{\delta}{\beta})^{2}+
(\frac{\delta}{\beta})\tilde{r}$ for some $a,b,\delta,\beta \in {\mathbb F}$ where $(a,b) \neq (0,0)$ and $\beta \not =0$. We let $\alpha =\beta/((\frac{b}{\beta})^{2}+
(\frac{b}{\beta})(\frac{a}{\beta})r+(\frac{a}{\beta})^{2}s)$. Consider the new
standard basis 
$$\begin{array}{ll}
 \tilde{x}_{1}=\frac{1}{\beta^{2}}[\alpha ar+\alpha b+\delta a)x_{1}+
  (\delta b+\alpha as)y_{1}], & \tilde{y}_{1}=\frac{1}{\beta}(ax_{1}+by_{1}) \\
 \tilde{x}_{2}=\frac{1}{\alpha\beta}(\alpha x_{2}+\delta x_{2}), &
    \tilde{y}_{2}=\beta y_{2}, \\
\tilde{x}_{3}=\frac{1}{\alpha}x_{3}, & 
\tilde{y}_{3}=\alpha y_{3}+\delta y_{2} \\
\tilde{x}_{4}=\frac{1}{\beta}[(\alpha ar+\alpha b +\delta a)x_{4}+
(\delta b+\alpha as)x_{5}], & \tilde{y}_{4}=\frac{1}{\beta^{2}}(by_{4}+\alpha
y_{5}), \\
\tilde{x}_{5}=ax_{4}+bx_{5}, & 
\tilde{y}_{5}=\frac{1}{\beta^{3}}[(\alpha ar+\alpha b+\delta a)y_{5}+
(\delta b+\alpha as)y_{4}].
\end{array}$$
Calculations show that $L$ has then presentation ${\mathcal P}_{10}^{(3,9)}(
\tilde{\gamma},\tilde{r},\tilde{s})$ with respect to the new standard
basis. \\ \\
It remains to see that the conditions are also necessary. Consider an algebra $L$ with presentation ${\mathcal P}_{10}^{(3,9)}(\gamma ,r,s)$ with respect to
some standard basis $x_{1},y_{1},\ldots ,x_{5},y_{5}$. Take some arbitrary new
standard basis $\tilde{x}_{1},\tilde{y}_{1},\ldots ,\tilde{x}_{5},
\tilde{y}_{5}$ such that $L$ has presentation 
${\mathcal P}_{10}^{(3,9)}(\tilde{\gamma},\tilde{r},\tilde{s})$ with respect to
the new basis where $\tilde{\gamma},\tilde{r}\not =0$ and where
$t^{2}+\tilde{r}t+\tilde{s}$ is irreducible. Then $\tilde{x}_{5}=ax_{4}+
bx_{5}$, $\tilde{y}_{2}=\beta y_{2}+u$ and $\tilde{y}_{3}=\alpha y_{3}+\delta y_{2}+v$
for some $u,v\in L^{3}$, $\alpha ,\beta ,\delta \in {\mathbb F}$ and
$(a,b) \in {\mathbb F} \times {\mathbb F} \setminus \{ (0,0)\}$ where $\alpha,\beta \neq 0$. The reader can convince himself that new basis that we
get satisfies the conditions stated. $\Box$ \\ \\
Before we give an example of an algebra of this form, we list some useful
properties of the groups $G(r,s)$ and $H(r)$. 
\begin{lemm} For any irreducible polynomials $t^{2}+rt+s$ and $t^{2}+
\tilde{r}t+\tilde{s}$, we have that \\ \\
(1) $H(\tilde{r})=(\tilde{r}/r)^{2}H(r)$. \\
(2) $G(\tilde{r},\tilde{s})=G(r,s)$ if $\tilde{s}-(\tilde{r}/r)^{2}s\in H(\tilde{r})$. 
\end{lemm}
{\bf Proof}\ \ Straightforward calculations. \\ \\
{\bf Example}. Let ${\mathbb F}$ be the finite field of order $2^{n}$. Let $r,s,
\tilde{r},\tilde{s}$ be as in the last lemma. Then $G(\tilde{r},\tilde{s})=
G(r,s)={\mathbb F}^{*}$ and thus $\tilde{r}/r\in G(r,s)$. Also
$[{\mathbb F}:H(\tilde{r})]=2$ and thus $\tilde{s}-(\tilde{r}/r)^{2}s\in
H(\tilde{r})$. It follows from the last proposition that the presentations
${\mathcal P}_{10}^{(3,9)}(\gamma ,r,s)$ and ${\mathcal P}_{10}^{(3,9)}
(\tilde{\gamma},\tilde{r},\tilde{s})$ describe the same algebra if and
only if $\tilde{\gamma}/\gamma\in ({\mathbb F}^{*})^{3}$. There are thus
either three algebras or one algebra according to whether $3$ divides 
$2^{n}-1$ or not. \\ \\
We end this section by giving a direct explanation why the relation 
     $$(\tilde{r},\tilde{s})\sim (r,s)\mbox{\ if\ }
       \frac{\tilde{r}}{r}\in G(r,s), \tilde{s}-(\frac{\tilde{r}}{r})^{2}s\in
      H(\tilde{r})$$
is an equivalence relation. \\ \\
First it is easy to see that $(r,s)\sim (r,s)$ as $1\in G(r,s)$ and $0\in H(r)$. Next
if $(\tilde{r},\tilde{s})\sim (r,s)$ then, as $G(r,s)=G(\tilde{r},\tilde{s})$
is a group, we have that $r/\tilde{r} \in G(\tilde{r},\tilde{s})$
and $s-(r/\tilde{r})^{2}\tilde{s}=(r/\tilde{r})^{2}(-\tilde{s}+
(\tilde{r}/r)^{2}s)\in (r/\tilde{r})^{2}H(\tilde{r})=H(r)$. This shows
that $\sim $ is symmetric. Finally suppose $(r^{*},s^{*})\sim (\tilde{r},\tilde{s})$ and $(\tilde{r},\tilde{s})\sim (r,s)$. Then we
have that $r^{*}/r=r^{*}/\tilde{r}\cdot
\tilde{r}/r \in G(r,s)$ and $s^{*}-(r^{*}/r)^{2}s=
[s^{*}-(r^{*}/\tilde{r})^{2}\tilde{s}]+[(r^{*}/\tilde{r})^{2}\tilde{s}-(r^{*}/r)^{2}s]=[s^{*}-(r^{*}/\tilde{r})^{2}\tilde{s}]+(r^{*}/\tilde{r})^{2}[\tilde{s}-(\tilde{r}/r)^{2}s]$ is
in $H(r^{*})+(r^{*}/\tilde{r})^{2}H(\tilde{r})=H(r^{*})$. Hence $\sim$ is also transitive and we have an equivalence relation.  

\bibliographystyle{plainurl}
\bibliography{biblio}
\end{document}